\documentclass[11pt, a4paper, twoside]{amsart}
\usepackage{amsmath}
\usepackage{amssymb}
\usepackage{amsfonts}
\usepackage{a4}
\usepackage[latin1]{inputenc}

\numberwithin{equation}{section}

\newtheorem{theorem}{Theorem}[section]

\theoremstyle{definition} 

\theoremstyle{plain} 

\newtheorem{corollary}[theorem]{Corollary}
\newtheorem{fact}[theorem]{Fact}
\newtheorem{conclusion}[theorem]{Conclusion}
\newtheorem{convention}[theorem]{Convention}
\newtheorem{claim}[theorem]{Claim}

\newtheorem{crucialfact}[theorem]{Crucial Fact}

\newtheorem*{maintheorem*}{Main Theorem} 
\newtheorem*{conjecture*}{Conjecture} 
\newtheorem*{claim*}{Claim} 
\newtheorem{definition}[theorem]{Definition}
\newtheorem{choice}[theorem]{Choice}
\newtheorem{hypothesis}[theorem]{Hypothesis}
\newtheorem{remark}[theorem]{Remark}
\newtheorem{notation}[theorem]{Notation}

\theoremstyle{remark}  
\newtheorem{remark*}{Remark}  

\newcommand{\nc}{\newcommand}
\nc{\nothing}[1]{}

\nothing{   
\nc{\dom}{{\rm dom}}
\nc{\card}{{\rm card}}
\nc{\lh}{{\rm lh}}
\nc{\lgg}{{\rm lg}}
\nc{\rge}{\mbox{\rm range}}
\nc{\cf}{{\rm cf}}
\nc{\nex}{\mbox{\rm next}}
\nc{\uhr}{\restriction}
\nc{\supt}{{\rm supt}}
\nc{\supp}{{\rm supp}}
\nc{\Lim}{{\rm Lim}}
\nc{\Leb}{{\rm Leb}}
\nc{\modd}{{\rm mod}}
\nc{\RO}{{\rm RO}}
\nc{\prob}{{\rm Prob}}
}

\nc{\On}{{\rm On}}

\nc{\nco}{\DeclareMathOperator}

\nco{\order}{o} 
\nco{\ppower}{pp} 
\nco{\pcf}{pcf} 
\nco{\tcf}{tcf} 
\nco{\tlim}{tlim} 
\nco{\limtext}{lim} 
\nco{\prodt}{{\textstyle \prod}} 
\nco{\symdiff}{\triangle}
\nco{\dom}{dom}
\nco{\card}{card}
\nco{\lh}{lh}
\nco{\lgg}{lg}
\nco{\rge}{range}
\nco{\otp}{otp}
\nco{\trunk}{tr}
\nco{\cf}{cf}
\nco{\nex}{next}
\nc{\uhr}{\restriction}
\nco{\supt}{supt}
\nco{\supp}{supp}
\nco{\Lim}{Lim}
\nco{\Leb}{Leb}
\nco{\modd}{mod}
\nco{\invariant}{inv}
\nco{\id}{id}
\nco{\RO}{RO}

\nc{\potom}{\ensuremath{{\cal P}(\omega)}}
\nc{\potinf}{\ensuremath{[\omega]^\omega}}
\nc{\pfin}{\ensuremath{{\cal P}(\omega)/{\rm fin}}}

\nc{\potfin}{\ensuremath{[\omega]^{<\omega}}}
\nc{\inn}{\ensuremath{{\omega^{\uparrow \omega}}}}
\nc{\hoch}{^{<\omega}}
\nc{\hocho}{^{\omega}}
\nc{\tree}[1]{{[} #1 {]}_0}
\nc{\tre}[2]{ {#1}_{#2}}
   
\nc{\prooff}[1]{{\bf Proof} of #1:}
\nc{\proofend}{\makebox{} \hfill ${\bf \square}$ \\}
\nc{\proofendof}[1]{\makebox{} \hfill $\boldmath{\square}_{\rm #1}$ \\}
\nc{\beq}{\begin{eqnarray*}}
\nc{\eeq}{\end{eqnarray*}}
\nc{\bde}{\begin{list}}
\nc{\ede}{\end{list}}

\newenvironment{myrules}
{\begin{list}{}
{
 \setlength{\leftmargin}{0.5in}
 \setlength{\labelwidth}{1cm}
 \setlength{\labelsep}{0.2in}
 \setlength{\parsep}{0.5ex plus 0.2ex minus 0.1 ex}
 \setlength{\itemsep}{0.3ex plus 0.2 ex minus 0ex}
}}{\end{list}}

{\end{list}}

\newcounter{subalph}
{\end{list}}

\newcommand{\greek}[1]{\ifthenelse{\value{#1}=1}{\mbox{$\alpha$}}%
  {\ifthenelse{\value{#1}=2}{\mbox{$\beta$}}{%
   \ifthenelse{\value{#1}=3}{\mbox{$\gamma$}}{%
   \ifthenelse{\value{#1}=4}{\mbox{$\delta$}}{%
   \ifthenelse{\value{#1}=5}{\mbox{$\varepsilon$}}{%
   \ifthenelse{\value{#1}=6}{\mbox{$\zeta$}}{%
   \ifthenelse{\value{#1}=7}{\mbox{$\eta$}}{%
   \ifthenelse{\value{#1}=8}{\mbox{$\theta$}}{%
   \ifthenelse{\value{#1}=9}{\mbox{$\iota$}}{%
   \ifthenelse{\value{#1}=10}{\mbox{$\kappa$}}{%
   \ifthenelse{\value{#1}=11}{\mbox{$\lambda$}}{%
   \ifthenelse{\value{#1}=12}{\mbox{$\mu$}}{%
   \ifthenelse{\value{#1}=13}{\mbox{$\nu$}}{%
   \ifthenelse{\value{#1}=14}{\mbox{$\xi$}}{%
   \ifthenelse{\value{#1}=15}{\mbox{$\rm o$}}{%
   \ifthenelse{\value{#1}=16}{\mbox{$\pi$}}{%
   \ifthenelse{\value{#1}=17}{\mbox{$\varrho$}}{%
   \ifthenelse{\value{#1}=18}{\mbox{$\sigma$}}{%
   \ifthenelse{\value{#1}=19}{\mbox{$\tau$}}{%
   \ifthenelse{\value{#1}=20}{\mbox{$\upsilon$}}{%
   \ifthenelse{\value{#1}=21}{\mbox{$\varphi$}}{%
   \ifthenelse{\value{#1}=22}{\mbox{$\chi$}}{%
   \ifthenelse{\value{#1}=23}{\mbox{$\psi$}}{\mbox{$\omega$}%
  }}}}}}}}}}}}}}}}}}}}}}}}

\newcounter{subgreek}
{\end{list}}

\newcounter{subarabic}
{\end{list}}

\newcounter{subroman}
{\end{list}}

\newcount\skewfactor

\def\mathunderaccent#1#2 {\let\theaccent#1\skewfactor#2
\mathpalette\putaccentunder}
\def\putaccentunder#1#2{\oalign{$#1#2$\crcr\hidewidth
\vbox to.2ex{\hbox{$#1\skew\skewfactor\theaccent{}$}\vss}\hidewidth}}
\def\name{\mathunderaccent\tilde-3 }
\def\Name{\mathunderaccent\widetilde-3 }

\nc{\nname}{\name}

\nc{\even}{\ensuremath{\rm Even}}
\nc{\odd}{\ensuremath{\rm Odd}}

\nc{\al}{$\alpha$\  }
\nc{\om}{\omega}
\nc{\omm}{\ensuremath{\omega_1}}
\nc{\ep}{\varepsilon}
\nc{\tk}{\tilde{K}}
\nc{\concat}{\,\hat{} \,}   
\nc{\force}{\Vdash}
\nc{\fb}{f_{\bar{M}}}
\nc{\such}{\, : \,}   

\nc{\meager}{\ensuremath{{\cal M}}}
\nc{\lebesgue}{\ensuremath{{\cal N}}}
\nc{\nulll}{\ensuremath{{\cal N}}}
\nc{\ksigma}{\ensuremath{{\bf K}_\sigma}}
\nc{\ideal}{\ensuremath{{\cal I}}}
\nc{\ga}{\ensuremath{\frak a}}
\nc{\AAA}{{\cal A}}   
\nc{\gc}{\ensuremath{\frak c}}
\nc{\gs}{\ensuremath{\frak s}}
\nc{\gh}{\ensuremath{\frak h}}
\nc{\gd}{\ensuremath{\frak d}}
\nc{\gb}{\ensuremath{\frak b}}
\nc{\gro}{\ensuremath{\frak g}}
\nc{\gu}{\ensuremath{\frak u}} 
\nc{\gr}{\ensuremath{\frak r}} 
\nc{\gt}{\ensuremath{\frak t}}
\nc{\fff}{\ensuremath{\frak f}}
\nc{\gm}{\ensuremath{\mathfrak{mcf}}}
\nc{\gge}{\ensuremath{\mathfrak e}}
\nc{\cfupro}{\ensuremath{\cf(\upro)}}
\nc{\cfvpro}{\ensuremath{\cf(\vpro)}}
\nc{\gp}{\ensuremath{\frak p}}
\nc{\gk}{\ensuremath{\frak k}}

\nc{\add}[1]{\mbox{\ensuremath{{\rm add}(#1)}}}
\nc{\cov}[1]{\mbox{\ensuremath{{\rm cov}(#1)}}}
\nc{\unif}[1]{\mbox{\ensuremath{{\rm unif}(#1)}}}
\nc{\cof}[1]{{\mbox{\ensuremath{\rm cof}(#1)}}}

\nc{\addd}[2]{\mbox{\ensuremath{{\rm add}^{#1}(#2)}}}   
\nc{\covv}[2]{\mbox{\ensuremath{{\rm cov}^{#1}(#2)}}}   
\nc{\uniff}[2]{\mbox{\ensuremath{{\rm unif}^{#1}(#2)}}} 
\nc{\coff}[2]{{\mbox{\ensuremath{\rm cof}^{#1}(#2)}}}

\nc{\cd}{Cicho\'n's Diagram}
\nc{\COF}{\mbox{\bf Cof}}

\nc{\MA}{\mbox{\rm MA}}
\nc{\PFA}{\mbox{\rm PFA}}
\nc{\OCA}{\mbox{\rm OCA}}
\nc{\GCH}{\mbox{\rm GCH}}
\nc{\CH}{\mbox{\rm CH}}
\nc{\zfc}{\mbox{\rm ZFC}}
\nc{\sch}{\mbox{\rm SCH}} 
\nc{\ZF}{\mbox{\rm ZF}}
\nc{\NCF}{\mbox{\rm NCF}} 
\nc{\FD}{\mbox{\rm FD}}   

\nc{\fourG}{\mbox{\rm 4G}} 
\nc{\fourI}{\mbox{\rm 4I}}   

\nc{\Borelhood}{Borel measurability} 
\nc{\Pieinseins}{\mbox{${\bf \Pi}^1_1$}}
\nc{\seinseins}{\mbox{${\bf\Sigma}^1_1$}}
\nc{\seinszwei}{\mbox{${\bf\Sigma}^1_2$}}
\nc{\seinsdrei}{\mbox{${\bf\Sigma}^1_3$}}
\nc{\Deleinszwei}{\mbox{${\bf\Delta}^1_2$}}

\nc{\up}{\ensuremath{{\cal U}\mbox{\ensuremath{\rm -prod}}\,\omega}}
\nc{\upp}{\ensuremath{{\cal U}'\mbox{\ensuremath{\rm -prod}}\,\omega}}
\nc{\upro}{\ensuremath{{\cal U}\mbox{\ensuremath{\rm -prod}}\,\om}}
\nc{\fupro}{\ensuremath{f({\cal U})\mbox{\ensuremath{\rm -prod}}\,\om}}
\nc{\vpro}{\ensuremath{{\cal V}\mbox{\ensuremath{\rm -prod}}\,\om}}
\nc{\fpro}{\ensuremath{{\cal F}\mbox{\ensuremath{\rm -prod}}\,\om}}

\nc{\cff}[1]{{\text{cf}\,(#1)}}           
\nc{\cu}{\ensuremath{\cal U}}             
\nc{\ai}{\ensuremath{\forall^\infty}}     
\nc{\ei}{\ensuremath{\exists^\infty}}     
\nc{\ww}{\ensuremath{\omega^\omega}}      

\nc{\gw}{groupwise dense}

\nc{\kk}{car\-dinal cha\-rac\-teris\-tic}
\nc{\joker}{\ast}
\nc{\gtc}{Galois-Tukey connection} 

\nc{\av}[1]{{\rm Av}_{#1}}
\nc{\eps}{\varepsilon}
\nc{\n}{{\bf n}}                 
\nc{\m}{{\bf m}}

\nc{\marginparr}[1]{}
\nc{\footnoteee}{} 
\nc{\footnotee}{}  

\newcommand{\cal}{\mathcal}

\nc{\divs}{{c_0 \setminus \ell^1}}
\nc{\divser}{(\divs, \leq^*)/\thickapproy}
\nc{\bfin}{\RO(\pfin \setminus\{0\},\subseteq^*)}
\nc{\bdivser}{\RO(\divser)}
\nc{\inc}{{\rm INC}}
\nc{\com}{{\rm COM}}
\nc{\thickapproy}{\makebox{}\!\!\thickapprox}
\nc{\approy}{\makebox{}\!\!\approx}
\nc{\lessi}{\leqslant}
\nc{\gessi}{\geqslant}
\nc{\interior}[1]{{\rm int}(#1)}
\nc{\closure}[1]{{\rm cl}(#1)}
\nc{\Vo}{Vojt\'a\v{s}}

\nc{\precedeseq}{\leq^*} 
\nc{\precedes}{\prec}
\nc{\stronger}{\leqslant_{\bf P}}
\nc{\underlline}[1]{\hat{#1}}
\nc{\PO}{{\bf P}}
\nc{\charak}{\text{ch}}
\nc{\needed}{needed\ }
\nc{\neededc}{needed}
\nc{\Needed}{Needed\ }
\nc{\wneeded}{weakly needed\ }
\nc{\Wneeded}{Weakly needed\ }
\nc{\wneededc}{weakly needed}

\begin{document}


\title{On needed reals}

\author{Heike Mildenberger and Saharon Shelah}
\thanks{The first author was supported by
a Minerva fellowship.}

\thanks{The second author's research 
was partially supported by the ``Israel Science
Foundation'', administered by the Israel Academy of Science and Humanities.
This is  the second author's work number 725}

\address{Heike Mildenberger,
Saharon Shelah,
Institute of Mathematics,
The Hebrew University of Jerusalem,
 Givat Ram,
91904 Jerusalem, Israel
}

\email{heike@math.huji.ac.il}
 
\email{shelah@math.huji.ac.il}

\begin{abstract}
Following Blass \cite{bl-ober},
we call a real $a$ ``\neededc'' for  a binary relation $R$ on the reals
if in every $R$-adequate set we find an
element from which $a$ is Turing computable. 
We show that every real needed for
${\bf Cof}({\lebesgue})$ is hyperarithmetic.
Replacing ``$R$-adequate'' by 
``$R$-adequate with minimal cardinality'' we get
related notion of being ``\wneededc''. We
show that is is consistent that the two notions do not coincide
for the reaping relation.
(They coincide in many models.)
 We show that not all hyperarithmetical reals are
\needed for the reaping relation.
This answers some questions asked by Blass at the Oberwolfach conference 
in December 1999 and in \cite{bl-ober}.
\end{abstract}



\maketitle

\tableofcontents

\nc{\LOC}{{\mathbb L}}

\setcounter{section}{-1}
\section{Introduction}

We consider some aspects of the following notions:

\begin{definition}\label{0.1}

\begin{myrules}

\item[(1)] (\Needed reals). Suppose that
we have a cardinal characteristic ${\mathfrak x}$ of the reals of the 
following form: There are (in most cases: Borel)
sets $A_-, A_+ \subseteq {\mathbb R}$ and there is
a (in most cases: Borel) relation $R \subseteq A_- \times A_+$ such that 
\[
{\mathfrak x}= ||R|| :=  \min \{ |Y| \such Y \subseteq A_+ \,\wedge\,
(\forall x \in A_-) (\exists y \in Y) R(x,y) \}.\]

We call $||R||$ the norm of $R$.
A set $Y \subseteq A_+$ is called $R$-adequate if $(\forall x \in \dom(R))
\; (\exists y \in Y) xRy$.
We say that $\eta \in {}^\omega 2$ is \needed
for $R$ if 
for every $R$-adequate set $Y$ there is some 
$y \in Y$ such that $\eta \leq_T y$.

If $A_+ \not\subseteq {\mathbb R}$ but can be mapped continuously
and injective and computably into $\mathbb R$ 
by a mapping $c$, called a coding, then we call the real $a$  \needed 
for $R$ and $c$ if
for any $R$-adequate set $Y \subseteq A_+$ there
is some $y \in Y$ such that $a \leq_T c(y)$.
In this situation, a real $a$ is called 
 \needed for $R$, if it is \needed for
$R$ and $c$ for any 
coding $c$.

\nc{\witnessing}{witnessing}

\item[(2)] (\Wneeded Reals). 
We call a real $a$ 
\wneeded for $R$ if for any $R$-adequate set $Y$
{\em of minimal cardinality} there is some
$y \in Y$ such that $a \leq_T y$. 

\end{myrules}
\end{definition}

Every \needed real is \wneededc. Sections~3 to~6 will
give some information on the reverse direction.

\section{\Needed reals for ${\COF}({\mathcal N})$}
\label{S1}

In this section we  answer
 Blass' question whether only hyperarithmetic reals
are \needed for the cofinality relation on the ideal of 
Lebesgue null sets affirmatively. 

\smallskip

In this section we work with two particular relations on the reals:
For functions $f,g  \colon \omega \to \omega$ we write $f
\leq^* g$ and say $g$ eventually dominates $f$ if
$(\exists n<\omega)(\forall k \geq n)(f(k) \leq g(k))$.
The dominating relation is
$${\bf D} = \{ (f,g) \such f,g \in {}^\omega \omega
\; \wedge \; f \leq^* g \},$$
and the cofinality relation for the ideal of
sets of Lebesgue measure zero is
$$\COF({\mathcal N}) = \{ (F,G) \such F,G 
\mbox{ are  $F_\sigma$-sets of Lebesgue measure 0 and } F \subseteq G \}.
$$
We write  ${\bf cof}({\mathcal N})$ for
$||\COF({\mathcal N})||$.

\medskip

Before stating our first theorem,  we review some
notation: For $s\in {}^{\omega>} 2 =
\{ t \such (\exists m \in \omega) t \colon m \to 2\}$,
we write $\lg(s)=\dom(s)$.
If $s \in {}^{\omega>} 2$ and  $ t \in 
{}^{\omega \geq} 2$,  we write 
$s \trianglelefteq t$ if $s = t \restriction \lg(s)$. 
Let  $s \triangleleft t$ denote that 
$s \trianglelefteq t$ and $s \neq t$.
A subset $T \subseteq {}^{\omega >}2$ is called a tree 
if it is
downward closed, i.e., if for all $t \in T$ for all 
$s \trianglelefteq t$,
we have that $s \in T$.
We let $\lim(T) = \{ f \in {}^\omega 2 \such 
(\forall n \in \omega) f \restriction n \in T \}$. 
An element $s \in T$ is a leaf if there is no 
$t \in T$ such that $s \triangleleft t$.
For a tree $T \subseteq {}^{<\omega}2$ and some $\rho \in
{}^{<\omega}2$ we set $T^{[\rho]} = \{ s \in T \such
s \trianglelefteq \rho \vee \rho \trianglelefteq s \}$.

$\Leb$ denotes the Lebesgue measure on the 
measurable subsets of ${}^\omega 2$, 
the product space of $\omega$ copies of the space $\{0,1\}$, where each
point has measure $\frac{1}{2}$.

\smallskip

We deal with the following forcings, where the
first is the ordinary Amoeba forcing. 
\begin{equation*}
\begin{split}
{\mathbb Q}  = &\biggl\{ T \such T \subseteq {}^{\omega >} 2,
T \mbox{ is a tree 
and  } \Leb(\lim(T)) > \frac{1}{2} \biggr\},\\
\underlline{\mathbb Q} = &
\left\{ T \in {\mathbb Q} \such \lim \left\langle \frac{|T \cap {}^n 2|}{2^n}
\such n \in \omega \right\rangle > \frac{1}{2}
\mbox{ and $T$ has no leaves} \right\},\
\intertext{ We set $ h_T(\rho) =
\Leb(\lim(T^{[\rho]}))$.}  
{\mathbb Q}^- = &
\left\{ T \in {\mathbb Q} \such \left( \forall
n \in \omega)( \rho \in {}^n 2 \cap T \rightarrow h_T(\rho) \cdot 2^{2^n} 
\in \omega \setminus \{0\}\right) \right\},\\
\underlline{\mathbb Q}^-  = & \underlline{\mathbb Q} \cap
{\mathbb Q}^-. 
 \end{split}
\end{equation*}

The partial order on ${\mathbb Q}$ and its variants is inclusion: 
subtrees are stronger ($\geq$, we follow the
Jerusalem convention) conditions.
It is easy to see that $\underlline{\mathbb Q}$, ${\mathbb Q}^-$
and  $\underlline{\mathbb Q}^-$ are
dense suborders of $\mathbb Q$.

\begin{theorem} 
\label{1.1}
Let $G$ be $\underlline{\mathbb Q}$-generic over $V$. Then in
$V[G]$ the following holds:
For every $\eta
\in {}^\omega 2\cap V$, if $\eta$ is recursive in the
generic tree $T=\bigcap G$, then $\eta$ is 
 \needed for domination.
\end{theorem}

\begin{conclusion}\label{1.2}
Since
 being \needed for domination is a
an absolute notion (see \cite{Jockusch, solovay:hyper} or \ref{4.1}),
also in $V$, every $\eta$ such that  $\eta
\in {}^\omega 2\cap V$  is recursive in $V[G]$ in the
generic tree $T=\bigcap G$, is \needed
 for domination.
\end{conclusion}

\noindent {\em Proof of 1.1.}
 For some $p \in \underlline{\mathbb Q}$, $\eta \in {}^\omega 2$, 
both in $V$, and Turing machine
$M$ (w.l.o.g.\ also in $V$) we have that

\begin{equation}\tag{$\ast$}
\label{computation}
p \Vdash \mbox{``$M$ computes $\eta$ from $\name{T}$''.} 
\end{equation}

Let $n(\ast) \in [1,\omega)$ and $p^* \in \underlline{\mathbb Q}^- $ 
be such that
$p \leq p^*$ and $\Leb(\lim(p^*)) = \frac{1}{2} + \frac{1}{n(\ast)}$.
Then, by the Lebesgue density theorem 
(3.10 in \cite{Oxtoby}),  we may 
choose $m(\ast)$ such that
for any $m \geq m(\ast)$,

\begin{equation*}
\begin{split}
\frac{1}{2} + \frac{1}{n(\ast)} \leq
\frac{|p^* \cap ({}^{m}2)|}{2^{m}} & \leq
\frac{1}{2} + \frac{1}{n(\ast)} + \frac{1}{2^{n(\ast) +7}}.
\end{split}
\end{equation*}
 
In order to derive from \eqref{computation}
some computation of $\eta$ relative to a suitable member
of a given ${\bf D}$-adequate set, we shall work with the following trees.
 
\begin{definition}\label {1.3}
For $r \in \underlline{\mathbb Q}$ and $\eps > 0$, if
$\Leb(\lim(r)) \geq \frac{1}{2} + \eps$  let
\begin{equation*}
\begin{split}
 T^\eps_{r,n} = 
\biggl\{ (q \cap {}^{n>} 2, h_q \restriction {}^{n>}2 ) 
\such & r \leq q \in \underlline{\mathbb Q}^-,\\
&
\Leb(\lim (q)) \geq \frac{1}{2} + \eps,\,
\forall m \frac{|q \cap {}^m2|}{2^m} \geq \frac{1}{2} + \eps \biggr\}.
\end{split}
\end{equation*}
We set $T^\eps_r = \bigcup\{ T^\eps_{r,n} \such n \in \omega \}$.
For $x \in T^\eps_{r,n}$ we write $x = (x(1),x(2))$.
We order $T^\eps_r$ by $\leq_T$: $(q \cap {}^{n>}2, h_q \restriction
{}^{n>}2) \leq_T (q'\cap{}^{n'>}2, h_{q'} \restriction{}^{n'>}2)$
iff $n \leq n'$ and 
$q \cap {}^{n>}2 =q' \cap {}^{n>}2$
and  $h_{q'} \restriction {}^{n>}2 = h_{q} \restriction {}^{n>}2$.
Equivalently, we may consider $t \in T^\eps_{r,n} $ as a function
$t \colon q \cap {}^{n>} 2 \to {\mathbb R}$, $t(\rho ) =
h_q(\rho)$.
We equip $T^\eps_r$ with the tree topology given by $\leq_T$, i.e.,
basic open sets in the topology are
$\{ t \in T^\eps_r \such t \geq_T t_0\}$, $t_0 \in T^\eps_r$.
\end{definition}

These trees exhibit the following properties:
\begin{myrules}
\item[$(\ast)_0$] $T^\eps_r$ is a tree with finite levels, the $n$th 
level being $T^\eps_{r,n}$.
\item[$(\ast)_1$]
If $\langle t_n \such n \in \omega \rangle $ is an $\omega$-branch of 
$T^\eps_r$ then $\Leb(\lim(\bigcup t_n(1)) \geq \frac{1}{2} + 
\eps$ and if $\eps > 0$ then $\bigcup t_n(1) \in 
{\mathbb Q}$.

\item[$(\ast)_2$] Moreover, we have if $r_1 \leq r_2$ in 
$\underlline{\mathbb Q}$ \and 
$\Leb(\lim(r_2)) - \frac{1}{2} \geq \eps$, then
$T^\eps_{r_2} \subseteq T^\eps_{r_1}  $.

\item[$(\ast)_3$] If $\Leb(\lim(r)) \geq \frac{1}{2} + \eps$, 
$p^\ast \leq r \in
{\mathbb Q}^-$ and $n \in \omega$ and $\langle t_\ell 
\such \ell \in \omega \rangle$ is an $\omega$-branch of $T^\eps_r$,
then for some $m \in \omega$,
there is $t^\ast \subseteq \dom(t_m)$ (here we regard $t$'s as
functions) such that 
\begin{myrules}
 
\item[(a)] $\sum \{t_m(\rho) \such \rho \in t^\ast \cap {}^m 2
\} > \frac{1}{2}$.

\item[(b)] If $M$ runs with input $n$ and oracle $f_{m,t^\ast}$ it
will give the value $\eta(n)$,
where $f_{m,t^\ast} \colon {}^{m \geq} 2 \to \{0,1\}$,
$f_{m,t^\ast} (\rho) =1 \Leftrightarrow (\exists \nu 
\in t^\ast) (\rho \trianglelefteq \nu)$.
\end{myrules}

\item[$(\ast)_4$] Let  $g^{\eps,\langle t_\ell \, : \, \ell \in \omega
\rangle}(n) $ be the first $m > n$ as in $(\ast)_3$.
For every $n, k \in \omega$ the sets
$$
S_{n,k} = \left\{ \bigcup_{\ell \in \omega} t_\ell \such \langle t_\ell
\such \ell \in \omega \rangle \mbox{ is a branch of } T^\eps_r \;
\wedge g^{\eps,\langle t_\ell \, : \, \ell \in \omega
\rangle}(n) \leq k \right\} 
$$
are open sets in the compact tree $T^\eps_r$, 
and $T^\eps_r = \bigcup_{k \in \omega} S_{n,k}$ is a
union of an increasing sequence $\langle S_{n,k}\such k \in \omega
\rangle$. Hence there is 
$K$, such that
$S_{n,K} = T^\eps_r$ and hence $K \geq 
g^{\eps,\langle t_\ell \, : \, \ell \in \omega
\rangle}(n)$ for all branches $\langle t_\ell \such \ell \in \omega \rangle$ 
of $T^\eps_r$. We let $g^\eps(n)$ be the minimal such $K$.

\end{myrules}

Now we specify the following items:
\begin{myrules}
\item[($\alpha$)] We take some
 $g \colon \omega \to \omega$ is such that $(\forall n)
g^\eps(n) \leq g(n)$. 
Our aim is to show that $\eta$ is recursive
in such a $g$.

\item[($\beta$)] 
 $\eps = \frac{1}{4 n(\ast)}$, and $\eps'= \frac{3}{4 n(\ast)}$. 
We choose some $p^\ast$ as above
and some $\underlline{\mathbb Q}$-generic filter $G$ such that $p^\ast 
\in G$. We fix an $\omega$ branch of $T^\eps_{p^\ast}$ such that
$t_{g(\ell)}$ determines $\eta(\ell)$ and the part of the oracle needed for it
in the sense of $(\ast)_3$ and $(\ast)_4$,
and $t_{g(\ell)}(1)$ is an initial segment of a condition in $G$.

\item[($\gamma$)] $p^{\ast\ast} = \{ \rho \such \rho \in p^\ast \cap 
{}^{m(\ast)} 2 \; \vee \;  (\rho \in
{}^{\omega>} 2 \setminus {}^{m(\ast)}2 \wedge
\rho\restriction m(\ast) \in p^\ast) \}$.
\end{myrules}

\nothing{
\item[$(\ast)_5$] if $p^{\ast\ast} \leq r \in \underlline{\mathbb Q}$,
$\Leb(r) \geq \frac{1}{2} + \frac{1}{2n(\ast)}$, then
$r \cap p^{\ast\ast} \in \underlline{\mathbb Q}$ (not in
${\mathbb Q}^-$ in general) and $\Leb(\lim(r \cap p^{\ast}))
\geq \frac{1}{2} + \frac{1}{4n(\ast)}$.
\end{myrules}}

The proof of the following claim will finish the proof of
Theorem~\ref{1.2}.

\begin{claim*}
 For every $n \in \omega$, $k \in \{0,1\}$, the following 
are equivalent:

\begin{myrules}
\item[(i)]  $\eta(n)=k$,

\item[(ii)] 
for some $t^1 \in T^{\eps'}_{p^{\ast\ast}, g(n)}$
(--- and this is recursive in $g$ ---)
 for every $t^0$ satisfying
$t^0 \subseteq t^1$ and $ t^0 \in T^\eps_{p^{\ast\ast},g(n)}$ there is
$t^2 \subseteq t^0$ such that 
$(\ast)_3$ (a) $+$ (b) holds with $t^\ast = \dom(t^2)$ and value
$\eta(n)= k$.

\nothing{
there is some $k \in \omega$ such that 
for $\eps = \frac{1}{k} $
if $t \in T^{2\eps}_{p^{\ast\ast}, g_\eps(n)} $ then for all
$t' \subseteq t$, $\sum \{t(2)(\rho)  \such 
\rho \in t' \cap {}^{g(n)}2 \}> \frac{1}{2}$,
and there is a run of $M$ with input $n$, oracle
$f_{g(n),t'}$ giving the result $k$.}

\end{myrules}
\end{claim*}

Proof:
(i) to (ii):
\nothing{
Take $r \geq p^\ast$ such that $r \in G$ and 
 $$ r \Vdash \mbox{``$M$ computes
$\eta(n)=k$''}.$$
By density, we may assume that $r \in \underlline{\mathbb Q}^-$.
We take an $\omega$-branch of $T^\eps_{p^\ast}$ as in $(\alpha)$.
Now we take $t^1 = r \cap t_{g(n)}\cap {}^{g(n)>}2$.
}
We assume (i). We take $t^1 = p^{\ast\ast}\restriction g(n)$.
If $t^0 \subseteq t^1$, $t^0 \in T^{\eps'}_{p^{\ast\ast},g(n)}$ is given,
we may take $t^2=t^0$. Since any branch containing $t^0$ and
stronger than $p^\ast$ forces $\eta(n) = k$, we have by the definition
of $g(n)$, that the part below $g(n)$ suffices for the computation.
So $t^2$ acts as desired.

\medskip
\noindent
(ii) to (i):

Assume that $\eta(n) = 1-k$. As we have ``(i) $\Rightarrow$ (ii)''
 for this situation, there is some 
$s^1 \in T^{\eps'}_{p^{\ast\ast}, g(n) }$ such that for every $s^0 \subseteq
s^1$ with $s^0 \in T^\eps_{p^{\ast\ast}, g(n)}$ there is $s^2 \subseteq s^0$
such that the analogues of $(\ast_3)$ (a) and (b) hold
with $\eta(n) = 1-k$. We have $t^1$ as in (ii) for $\eta(n)= k$.
There are $q_0,q_1$ witnessing
$t^1, s^1 \in T^{\eps'}_{p^{\ast\ast}, g(n)}$.

Subclaim 1: $q_0, q_1$ are compatible in the Amoeba forcing.
\\
Proof of the claim: Both satisfy:

\begin{equation*} 
\begin{split} 
&\lim(p^{\ast\ast}) \supseteq \lim(q_\ell),\\
&\frac{1}{2} + \frac{1}{n(\ast)}  \leq
\Leb(\lim(p^{\ast\ast})) \leq \frac{|p^\ast \cap {}^{m(\ast)} 2|}{2^{m(\ast)}}
\leq \frac{1}{2} + \frac{1}{n(\ast)} + \frac{1}{2^{n(\ast)+7}},\\
&\Leb(\lim(q_\ell)) \geq \frac{1}{2} + \eps'.
\end{split}
\end{equation*}

We show that $\Leb(\lim(q_0) \cap \lim(q_1)) > \frac{1}{2}$:

We have that 
\begin{equation*}\begin{split}
\Leb(\lim(p^{\ast\ast}) \setminus & (\lim(q_0) \cap \lim(q_1)))\\
&\leq 
\Leb(\lim(p^{\ast\ast}) \setminus (\lim(q_0))) + 
\Leb(\lim(p^{\ast\ast}) \setminus  \lim(q_1))) 
\\
&\leq  2 \cdot \left(\frac{1}{4 n(\ast)} + \frac{1}{2^{n(\ast) + 7 }}
\right)
\\
& = \frac{1}{2n(\ast)} + \frac{1}{2^{n(\ast) + 6}},
\end{split}
\end{equation*}
hence
\begin{equation*}
\begin{split}
\Leb(\lim(q_0) \cap \lim(q_1))
\geq & 
\Leb(\lim(p^{\ast\ast})) -
\Leb(\lim(p^{\ast\ast}) \setminus (\lim(q_0)\cap \lim(q_1))) \\
\geq & \frac{1}{2} +  \frac{1}{ n(\ast)} 
- \frac{1}{2n(\ast)} - \frac{1}{2^{n(\ast) + 6}} > \frac{1}{2}.
\end{split}
\end{equation*}
So the sublcaim is proved.

But: $q_0$ and $q_1$ cannot be compatible in the Amoeba forcing.
By the choice of $p^{\ast}$ we have that
$$
p^\ast \Vdash \mbox{``} \eta \mbox{ is computed by $M$
using the oracle $\Name{T}$.''}
$$
We have that $q_\ell \geq p^{\ast}$ and that $q \geq q_\ell$.
But then we find $t^2_\ell \subseteq 
p^\ast \cap 2^{g(n)}$ such that 
\begin{myrules}
\item[(a)] $\sum_{x \in t^2_\ell} h_{p^\ast}(x) > \frac{1}{2}$, and
\item[(b)] if $M$ runs on the input $n$ and the oracle 
$t^2_\ell$ it will give the result $\eta(n)$  for $\ell =0$ and
$1-\eta(n) $ for $\ell =1$.
\end{myrules}
Since $\eta \in V$, there cannot be two different computations,
depending on two different $\name{T}[G] \cap {}^{g(n)} 2$.
Hence the assumption that $q_0$ and $q_1$ with the above properties 
both exist leads to a contradiction, and the Claim and 
Theorem~\ref{1.2} are proved.

\begin{theorem}\label{1.4}
Every \needed real  for ${\COF}({\mathcal N})$ is \needed for
the dominating relation.
\end{theorem}

\proof: Let $\{ A_i \such i < \kappa\}$ be 
a ${\COF}({\mathcal N})$-adequate set, such that 
each $A_i$ is a Borel set. Let $\eta
\in {}^\omega 2$.

For each $i$ choose a countable elementary submodel $N_i$ 
of $({\mathcal H}(\beth_3),\in)$ to which $\eta$ and $A_i$ belong. 
We let $G_i$ be  a 
subset of ${\mathbb Q}^{N_i}$ that is generic over $N_i$ and let $T_i =
\name{T}[G_i]$. Now let $A_i^\ast$ be
\begin{multline*} A_i^\ast =
\{\rho \in {}^\omega 2 \such \mbox{ no $\rho' \in {}^\omega 2$ 
which is almost equal to $\rho $ }\\
\mbox{(i.e. $\rho(n) = \rho'(n)$ for every large enough $n$) belongs 
to $T_i$}\}
\end{multline*}

$A_i^\ast$ is a null set: We have
$A_i^\ast = \bigcap_{n \in \omega } (\{ \rho' \such 
(\exists \rho \in T_i) \; (\rho'\restriction[n,\omega)
= \rho \restriction[n, \omega)) \})^c$. Furthermore we have that 
$\lim_{n \to \infty} 
\Leb(\{ \rho' \such (\exists \rho \in T) (\rho' \restriction [n,\omega)
= \rho \restriction [n,\omega)) \}) = 1$, because for a given
$\eps$, by the Lebesgue density
Theorem (3.20 in \cite{Oxtoby}) there is some $n_0$ such that for 
$n \geq n_0$ we have for all
$s \in T\cap {}^n 2$ that $\Leb(T \cap [s]) \cdot 2^n > 1-\eps$
and hence $
\Leb(\{ \rho' \such \exists \rho \in T \rho' \restriction [n,\omega)
= \rho \restriction [n,\omega) \}) > 1-\eps$.

By genericity of $T_i$ and because $A_i \in N_i$ and
because $A_i $ is a nullset in $N_i$ we have that $A_i \subseteq (T_i)^c$.
The same argument shows that for all $s \in 
{}^{\omega >}2$ we have that $\{ s \concat f \such 
\exists s' \; (|s'| = |s| \,\wedge\, s' \concat f \in A_i) \}$
is a subset of  $(\lim(T_i))^c$. Hence we have that
$A_i \subseteq A_i^\ast$. Therefore also $\{A^\ast_i \such
i < \kappa\}$ is a ${\COF}({\mathcal N})$-adequate set.
 If $\eta$ is recursive in $A_i^\ast$
(more precise: in one one of $A_i^\ast$'s simple codings) 
it is also recursive in $T_i$ and hence by Theorem~\ref{1.2}
\needed for dominating.
\proofend

\smallskip

\begin{fact}\label{1.5}
We use the result of Jockusch and Solovay
every real that is \needed  for the  dominating relation
is hyperarithmetic (Solovay \cite{solovay:hyper}) and this is 
optimal (Jockusch, \cite{Jockusch}): every hyperarithmetic real
is \needed for the dominating relation.
\end{fact}

Blass \cite[Theorem 6, Corollary 8]{bl-ober} showed that every real that is 
needed for ${\bf D}$ is also needed for ${\COF}({\mathcal N})$
and hence that all hyperarithmetic reals are needed for 
 ${\COF}({\mathcal N})$. So this gives the other inclusion 
in the following corollary:

\begin{corollary}\label{1.6}
Exactly the hyperaritmethic reals are needed for the 
${\COF}({\mathcal N})$-relation.
\end{corollary}

\nothing{
\medskip
We can axiomatize the proof.
\begin{theorem}\label{3}
Assume that for every $y \in \rge(R)$ and for every $\eta$ that is
not needed for dominating there are ${\mathbb Q}$ and $\name{Y}$ such that
\begin{myrules}
\item[(i)] ${\mathbb Q}$ is a forcing notion,
\item[(ii)] $\mathbb Q$ is ccc,
\item[(iii)] $\name{Y}$ is a $\mathbb Q$-name of a member of $\rge(R)$ 
such that $$ \forall x (x R y \rightarrow x R \name{Y})$$
(for countably many names)
\item[(iv)] for every $p \in {\mathbb Q}$ there is a recursive in
$g$  subtree $T^\ast$ of
${}^{\omega} \omega$ and there is a name $\name{T}$ such that
\begin{myrules}
\item[(a)] for every $\nu \in T^\ast$ 
some $q$, $p \leq q \in {\mathbb Q}$ and $q$ forces that 
$\nu$ is included in $\name{T}$,
\item[(b)] $p \Vdash \name{T} \in \lim T^\ast$,
\item[(b)] if $\eta$ is recursive in 
$\name{Y}$ then it is recursive in $\name{T}$.
\end{myrules}
\end{myrules}
\end{theorem}
\proof
This suffices because we can compute $\eta$ from $g$ by 
$\eta(n) = k$ iff for some $\nu \in T^\ast$ if we let $M$ run
with a machine according to item (c)
 with $\nu$ as an oracle the run finishes and gives result $k$.
}

\section{Needed reals for the slalom relation and a
general scheme}\label{S2}

In this section we deal with a forcing $\LOC$ which 
is closely related to the localization forcing from
\cite[page 106]{BJ}. Theorem~\ref{2.3} is analogous to Theorem~\ref{1.1},
but for the forcing ${\mathbb L}$. Theorem~\ref{2.10} is analogous to
Theorem~\ref{1.4}, but the translation mechanism in the proof is different.

In the second part of the section, we collect sufficient conditions and give a
general scheme for the proofs of ``being computable in the generic and
being in $V$ implies being hyperarithmetic'' and of ``every real
\needed for $R$ is $\Delta_1^1$.''

\begin{definition}\label{2.1}
\begin{equation*}
\begin{split}
\LOC = & \{ p \such p=(n,\bar{u})=(n^p,\bar{u}^p), \bar{u}=
\langle u_\ell \such \ell \in \omega \rangle, u_\ell
\in [\omega]^{\leq \ell}, 
\\
& h(p) := \limsup
\langle |u_\ell| \such \ell \in \omega \rangle < \omega
\mbox{ is well-defined}\},\\
p \leq q & \leftrightarrow \left(
\bigwedge_{\ell \in \omega } u_\ell^p \subseteq 
u_\ell^q \wedge \bar{u}^q \restriction n^p = \bar{u}^p \restriction 
n^p\right).
\end{split}
\end{equation*}
The generic is considered as a characteristic function
$\rho$ with domain $\omega \times \omega$ such that
$\rho(n,m) =1 \leftrightarrow (\exists p \in G)
(m \in u^p_n) $.

\end{definition}

\begin{notation}\label{2.2}
An $m$-oracle is a  function from $m \times m$ to $\{0,1\}$. If
$\bar{u} = \langle u_\ell \such \ell < m \rangle$, $u_\ell
\in [\omega]^{<\ell}$ the $\bar{u}$-oracle 
$\rho_{\bar{u}}\in {}^{m\times m} 2$ is defined by
$\rho_{\bar{u}}(n_1,n_2) =1 \leftrightarrow 
n_2 \in u_{n_1}$. We allow that $(\exists \ell<m)
\; \max(u_\ell) > \lg(\bar{u})=m$.
\end{notation}

\begin{theorem}\label{2.3}
 Assume that $M$ is a Turing machine and that $\eta \in {}^\omega 2$.
Let $\name{G}$ be a name 
for an $\LOC$-generic element.
Suppose that $p^\ast \in \LOC$ and that 
$$p^\ast \Vdash_{\LOC} M \mbox{ computes $\eta$ from  } \name{G}.
$$
Then $\eta$ is hyperarithmetic.
\end{theorem}

\proof
Let $n^\ast = n^{p^\ast}$ and
 $\bar{u}^\ast = \bar{u}^{p^\ast} \restriction n^{\ast}$,
and $h^\ast = h(p^\ast)$.
By a density argument we may 
 assume that $ n^\ast > 4 h^\ast \wedge (\forall \ell)(\ell 
\geq n^\ast \rightarrow |u_\ell^\ast|
\leq h^\ast)$.

We let 
\begin{equation*}
\begin{split}
T= T_{\bar{u}^\ast} =& \{ \bar{u} \such n^\ast \leq m < \omega, 
\bar{u} = \langle u_\ell \such \ell \in m \rangle, u_\ell
\in[\omega]^{\leq \ell},\\
&
\bar{u} \restriction n^\ast =
\bar{u}^\ast \wedge \ell \geq n^\ast \rightarrow |u_\ell| =
h^\ast \}.
\end{split}
\end{equation*}
We order $T$ by the initial segment relation $\trianglelefteq$.
The set of all infinite branches of $T$ is
$[T] = \{ \bar{u} \such \forall n \bar{u} \restriction n \in T \}$.

If $\bar{u} \in T_{\bar{u}^\ast}$  we let 
\begin{equation*}\begin{split}
\Xi_{\bar{u}} = & \Xi_{\bar{u},\bar{u}^\ast} =
\{ \rho_{\bar{v}} \such
\lg(\bar{u})= \lg(\bar{v}), \bar{v} \restriction n^\ast =
\bar{u}^\ast,\\
&
 (n^\ast \leq \ell< \lg(\bar{u})
\rightarrow [0,\lg(\bar{u})) \cap u_\ell \subseteq v_\ell)
\}.
\end{split}
\end{equation*}

\begin{fact}\label{2.4}
For every $j< \omega$, $\bar{u} = \langle  u_\ell \such \ell \in \omega 
\rangle \in [T_{\bar{u}^\ast}]$, such that for each $\ell$,
and $\lg(u_\ell) = \ell$, 
there are
$m \in [n^\ast,\omega) $, $\bar{v} \in \Xi_{\bar{u} \restriction m}
\cap \Xi_{\bar{u}^\ast}$ such that
\begin{equation}
\tag{$\ast$}
\begin{split}
&\mbox{ with $\rho_{\bar{v}}$ as an oracle on domain
$\lg(\bar{v}) \times \lg(\bar{v})$, $M$ finishes its run}
\\
&
\mbox{ and gives the result $\eta(j)$}
\end{split}
\end{equation}
\end{fact}

\proof
The conditions 
$(n^\ast,\bar{u})$ and $p^\ast= (n^{p^\ast}, \bar{u}^{p^\ast})$ 
are compatible: $(n^\ast,\bar{v})=(n^\ast, \langle u_\ell \cup
u^{p^\ast}_\ell \such \ell \in \omega \rangle ) \in \LOC$
is stronger or equal to both of them (here we use 
$n^\ast > 4 h^\ast$) and in 
$\Xi_{\bar{u}\restriction m}$ for all $m$. We take a generic to
which $(n^\ast,\bar{v})$ belongs.
Consider the run of $M$, it uses only
$\bar{v} \cap(m \times m)$ for  $m$ large enough.
\proofend

\begin{fact}\label{2.5}
For every $j< \omega$ 
there is $m_j \in (n^\ast,\omega)$ such that 
such that for every $
\langle u_\ell \such \ell \in \omega \rangle \in [T_{\bar{u}^\ast}]$,
there is $\rho_{\bar{v}}
 \in \Xi_{\bar{u} \restriction {m_j}}
\cap \Xi_{\bar{u}^\ast}$ such that
($\ast$) holds.
\end{fact}

\proof By the previous lemma and by K\"onig's lemma.
All the  levels of $T_{\bar{u}^\ast}$ are finite.   Note that
$\Xi_{\bar{u}}$ depends only on $\langle u_\ell \cap
\lg(\bar{u}) \such \ell 
< \lg(\bar{u}) \rangle $.
\proofend

\begin{definition}\label{2.6}
$g_{M,\bar{u}^\ast} \in {}^\omega \omega $ is defined by
$$g_{M,\bar{u}^\ast}(j) = \mbox{$\min\{ m_j \such m_j $ in as in the 
Fact~\ref{2.5}}\}.$$
\end{definition}

\begin{claim}\label{2.7}
For every $j \in \omega, k < 2$ and $m \geq g_{M,\bar{u}^\ast}(j)$ the
following are equivalent:
\begin{myrules}
\item[(i)]
$\eta(j) = k$,
\item[(ii)]  
for some $\bar{u} =\langle
u_\ell  \such \ell <m \rangle$ and $h^\ast$, 
such that $(\ell \in [n^\ast,m)
\rightarrow u_\ell \in [m]^{\leq h^\ast})$, $\bar{u}
\restriction n^\ast = \bar{u}^\ast$ for every 
$\bar{u}' =\langle
u'_\ell  \such \ell <m \rangle$ such that $\ell \in [n^\ast,m)
\rightarrow  u'_\ell \in [m]^{\leq h^\ast}$, $\bar{u}'
\restriction n^\ast = \bar{u}^\ast$ 
there is $\bar{v} \in \Xi_{\bar{u}} \cap \Xi_{\bar{u}'}
\subseteq {}^{m \times m}2 $ such that
$M$ running with oracle $\rho_{\bar{v}}$ and input $j$ 
finishes its run and gives the
result $k$.
\end{myrules}
\end{claim}
\proof:
(i) $\Rightarrow$ (ii): By the previous fact, $\bar{u}^{p^\ast} 
\restriction m$  is as required.
(ii) $\Rightarrow$ (i):
Let $\bar{u}$ be as guaranteed in (ii). It is said there ``for every 
$\bar{u}'$'' so in particular for $\bar{u}'= \bar{u}^{p^\ast} 
\restriction m$, 
there is $\rho \in \Xi_{\bar{u}} \cap \Xi_{\bar{u}'}$ as there. 
Now we can find a condition $q \in \LOC$ such that 
$n^q=m > n^\ast$, $\bar{u}^q \restriction n^\ast = \bar{u}^\ast$,
$n^\ast \leq \ell < m \Rightarrow u^q_\ell = u^{p^\ast}_\ell
\cup v_\ell = u'_\ell \cup v_\ell$,
$\ell \geq m \rightarrow u^q_\ell = u^{p^\ast}_\ell$.
So
\begin{myrules}
\item[($\alpha$)] $p^\ast \leq q$ and $q \Vdash \name{G} \restriction
{}^{m \times m} 2 = \rho_{\bar{v}}$, hence
\item[($\beta$)]
 $q \Vdash$ ``$M$ running with the oracle $\name{G}$
and input $j$ gives the result $k$'', and recall
\item[($\gamma$)] $p^\ast \Vdash $ ``$M$ computes $\eta$''.
\end{myrules}
By ($\alpha) + (\beta) + (\gamma)$ we get that $\eta(j)=k$ is as required.
\proofend

\begin{conclusion}\label{2.8}
Assume that $\eta \in {}^\omega \omega$, $\eta \in N$, $G$ is $\LOC$-generic
over  $N$ and  that $\name{\rho}[G] = \rho$ and $N[G] 
\models \mbox{``}\eta \leq_T \rho$''.
Then $\eta$ is hyperarithmetic.
\end{conclusion}

\proof Analogous to the proof of \ref{1.2} for
$N$ instead of $V$. We use \ref{2.7}.

\begin{definition}\label{2.9}
$S \in {}^\omega({}^{\omega >}[\omega])$ is called a slalom 
iff for all $n$, $|S(n)| \leq n$.
\end{definition}

\begin{theorem}\label{2.10}
Exactly the hyperarithmetic reals are \needed for
the slalom relation
\[{\bf SL}=
\{ (f, S) \such f \in {}^\omega \omega \; \wedge \;
 S \mbox{ is a slalom and }
(\forall n \in \omega) (f(n) \in S(n)) \}.\]
\end{theorem}

\proof
First show that only hyperarithmetic reals are
\needed for ${\bf SL}$:
Let $\{ S_i \such i < ||SK||\}$ be an $SL$-adequate set. 
Let $\eta \in {}^\omega 2$.
We take $N_i \prec (H(\beth_3),\in)$ such that $ \eta, S_i \in N_i$.
Then we let $G_i$ be ${\mathbb L}$-generic over $N_i$. Now we set
$S_i^\ast = \{ \rho \such (\exists \rho' \in G_i) 
\rho' =^\ast \rho \}$. Then we
have that $S_i \subseteq S_i^\ast$, $S_i^\ast$ is the union of $\omega$ 
slaloms, each of them computable from $G_i$, and the members of
all the unions form an ${\bf SL}$-adequate set.

All hyperarithmetic reals are \needed for ${\bf SL}$, because all of then
are \needed for ${\bf D}$. 
Suppose that $\{ \langle S^\alpha_i \such i \in \omega \rangle
\such \alpha \in ||{\bf SL}|| \}$ is ${\bf SL}$-adequate and that $\eta
\in {}^\omega 2$ is hyperarithmetic. Then
$\{ \langle \max{S_i} \such i \in \omega \rangle \such 
\alpha \in ||{\bf SL}||\}$ is  ${\bf D}$-adequate and hence there is some 
element $f$ in it from which $\eta$ is computable. But then of 
course $\eta$ is also computable in any slalom where $f$ stems from. 
\proofend

\medskip

From our two examples $(\mathbb{Q},{\bf Cof}({\lebesgue}))$ and
 $({\mathbb L}, {\bf SL})$ we collect
the following scheme:

\begin{theorem} \label{2.11}
Assume that 
\begin{myrules}
\item[(a)] $T \subseteq H(\aleph_0)$ is recursive, 
$T$ is a tree with $\omega$ 
levels and each level is finite, each $v \in T$ is a finite function from 
$H(\aleph_0)$ to $H(\aleph_0)$.

\item[(b)] $Q$ is a forcing notion, and $\rho_n$, $n \in \omega$,
are $Q$-names, and
$$\Vdash_Q (\forall n \in \omega)
\; (\rho_n \in \lim(T))
\; \wedge \; (\forall x \in \rge(R)) \bigvee_{n \in \omega}
\forall y (yRx \rightarrow yR \rho_n).$$

\item[(c)] For each $n \in \omega$ we have: For a dense set of $p_0 \in Q$ 
there is some $p \geq p_0$ such that
the following conditions are fulfilled:
\begin{myrules}

\item[($\alpha$)]
Let $T_{n,p} =
\{ \nu \in T \such 
p \Vdash \nu \subseteq \rho_n \}$. This is a subtree of $T$.

\item[($\beta$)]
 Let $S^\ast_{n,p} = \Bigl\{ t \such 
$ for some subtree $T'$ of $T_{n,p}$ 
and some $k$,
 $t=
\{ \nu \in T' \such $ level$_{T_{n,p}}(\nu) \leq k \}$, 
and no maximal node of $t$ has level $k \Bigr\}$, 
and order $S^\ast_{n,p}$ naturally.

\item[$(\gamma)$] 
$S_{n,p}$ is a recursive 
subtree of $S_{n,p}^\ast$ such that 
\begin{myrules}
\item[(i)] $T_{n,p}$ is an $\omega$-branch of $S_{n,p}$,
\item[(ii)] for every branch $\bar{t} =
\langle t_\ell \such \ell \in \omega \rangle$ of $S_{n,p}$ 
there is $q \in Q$
such that $q$ is compatible with $p$ and $T_{n,q}
= \bigcup_{\ell \in \omega} t_\ell$.
\end{myrules}
\end{myrules}
\item[(d)] $\eta \in {}^\omega 2$ or ${}^\omega \omega$
\end{myrules}

Then we have for every $n \in \omega$:
 if $\Vdash_Q$ ``$\eta$ is recursive in $\name{\rho_n}$'' 
then $\eta$ is hyperarithmetic.
\end{theorem}

\proof

So for some $p^\ast$ as in (c) and Turing machine $M$
$$p^\ast \Vdash_Q \mbox{``$M$ computes $\eta$
from $\name{\rho_n}$''}.
$$
Let $S_{n,p^\ast}$ and $S^\ast_{n,p^\ast}$ be as in clause (c). Now
we prove some intermediate facts, and the proof of \ref{2.11}
will be finished with \ref{2.15}.

\begin{fact}\label{2.12}
For every $\omega$-branch $\langle t_k \such k \in \omega \rangle$
of $S_{n,p^\ast}$ and $j \in \omega$ for some (= every) large enough 
$k \in \omega$ for some
$\nu \in t_k \cap {\text{level}}_k(T_{n,p^\ast})$
 if $M$ runs on input $j$ and oracle $\nu$ it finishes 
(so we do not ask oracle questions outside the domain) and gives the 
result $\eta(j)= k$.
\end{fact}

\proof
There is $q$ compatible with $p^\ast$ such that $T_{n,q}
\subseteq \bigcup_{n \in \omega} t_n$ . Let $r \geq p^\ast, q$, and let
$G \subseteq Q$ be generic with $r\in G$, so $p^\ast \in G$. If $M$ runs with
$\name{\rho_n}[G] $ it gives $\eta(j)$, so for some $\nu \in T$, $\nu 
\subseteq \name{\rho_n}[G]$. And $M$ can use as an oracle only $\nu$, but as 
$q \in G$, $\nu \in T_{n,q} \subseteq \bigcup_{\ell \in \omega} t_\ell$.
Of course any $\nu'$, $\nu \subseteq \nu' \in T_{n,p^\ast}$ can serve.

\begin{fact}\label{koenig}
For $j \in \omega$, for every large enough $m$, for every $t \in
{\text level}_m(S_n)$ there is $\nu \in t \cap 
{\text level}_m(T_{n,p^\ast})$ 
such that if $M$ runs with $\nu$ as an oracle then it computes $\eta(j)$.
\end{fact}
\proof By the previous fact and K\"onig's lemma.

\begin{definition}\label{2.14}
We define $g_{p^\ast} \in {}^\omega \omega$ by $g_{p^\ast}(j) =
\min\{m \such m \mbox{ as in } \ref{koenig}\}$.
\end{definition}

\begin{crucialfact}\label{2.15}
 For $j,n \in \omega$, $k \in 2$, the following are
 equivalent for any $m \geq g_{p^\ast}(j)$:
\begin{myrules}
\item[(i)] $\eta(j)=k$.
\item[(ii)] there is $t^1 \in {\text level}_m(S_{n,p^\ast})$ 
such that for every
$t^2 \in {\text level}_m(S_{n,p^\ast}) $ there is 
$\nu \in t^1 \cap t^2$ such that if we let run 
$M$ with input $j$ and oracle $\nu$ then the run finishes 
and there are no questions to the oracle that do not have an answer,
and it gives answer $k$.
\end{myrules}
\end{crucialfact}
\proof Analogous to \ref{2.7} \proofendof{2.11}
\begin{remark}\label{2.16}
1. Usually, $S_{n,p^\ast}$ is not so dependent on $p^\ast$, rather we have
that $Q= 
\bigcup_{k \in \omega} Q_k$, and for all $k \in \omega$ we have 
$S_{n,p^\ast}$ 
as above 
being the same for each $p^\ast \in Q_k$.

2. Actually we use in (c)($\gamma$)(i) only $T_{n,q} =
\bigcup_{k \in \omega} t_k$.
But we use $T_{n,p^\ast} = \bigcup_{k \in \omega} t_k$ for some 
$\omega$-branch. 
\end{remark}

\begin{theorem}\label{2.17}
A sufficient condition for ``every real needed for $R$ is $\Delta_1^1$''
is:
For some forcing notion $Q$ and some $
Q$-names $\name{\rho_n}$, $n \in \omega$, we have
\begin{myrules}
\item[(a)] $\Vdash_Q  \mbox{``}\name{\rho_n} 
\in [T], \name{\rho_n} \in \rge(R)$''
\item[(b)] $\Vdash_q $``for every $x \in \dom(R)$ for some $n$,
 $ xR \name{\rho_n}  $''
\item[(c)] for each $n$:
$Q$, $T$ and $\name{\rho_n}$ satisfy the conditions in \ref{2.10} or 
just its conclusion.
\end{myrules}
\end{theorem}
\proof Like the first half of the proof of Theorem~\ref{2.10}.

\section{\Wneeded reals for the reaping relation}\label{S3}

In this section we show that it is consistent that all hyperarithmetic reals
 are \wneeded for the reaping relation.
In Section 5 we shall prove in $\zfc$ that not all
hypearithmetic real are are \needed for the reaping relation, answering
another question from Blass' work \cite{bl-ober}.
In a model of $\CH$, the notions ``\needed real'' and ``\wneeded real'' 
coincide, and thus in such a model 
not all hyperarithmetic reals are  
\wneeded for the reaping relation. 
The model of this section, together with the result from Section 5,
gives an example for the 
fact that in contrast to the notion of ``being \neededc'', 
the notion of ``being
\wneededc'' is not absolute.

\begin{definition}\label{3.1}
 The relation 
$${\bf R} = \{ (f,X) \such f \in {}^\omega 2, X \in {}^\omega[\omega]
\wedge f \restriction X \mbox{ is constant}\}
$$ 
is called the  reaping or the refining or the
unsplitting relation. We say ``$X$ refines $f$'' if $f \restriction X$
is constant. We say ``${\mathcal R}$ refines $f$'' if there is some $X 
\in {\mathcal R}$ that refines $f$. Finally we say
``${\mathcal R}$ refines $F$'' if for every $f \in F$ we have that
${\mathcal R}$ refines $f$.

The norm of this relation is called $\gr$, the reaping number
or the refining number or the 
unsplitting number.
\end{definition}

\begin{definition}\label{3.2}
Let $g \in {}^\omega \omega$ be strictly increasing and $g(n) > n$.

\begin{myrules}
\item[(1)] We say $A \in [\omega]^\omega$ is $g$-slow if $(\exists^\infty n) 
|A \cap g(n) | \geq n$.
\item[(2)] 
\begin{multline*}
{\mathcal F}_g = \{ f \such \dom(f) \in [\omega]^\omega,
\mbox{ for $i \in \dom(f)$ we have that } f(i)=(f^1(i),f^2(i))\\
\mbox{ and } 
 f^2(i) \in [g(f^1(i))]^{\geq f^1(i)} \mbox{ and } 
\limsup \langle f^1(i) \such i \in \dom(f) \rangle
= \omega \}.
\end{multline*}
\item[(3)]
 We say that $\bar{A}$ is $(g,\kappa)$-o.k.\ if 
\begin{myrules}
\item[(a)] $\bar{A} = \langle A_i \such i < \kappa \rangle$, and 
\item[(b)] $A_i \in [\omega]^\omega$,
\item[(c)] if $k < \omega$, $f_0, \dots , f_{k-1} \in {\mathcal F}_g$,
$\bigcap_{\ell \in \omega} \dom(f_\ell) = B \in [\omega]^\omega$
and $\limsup\langle \min\{f^1_\ell(i) \such \ell \in k \}
\such i \in B \rangle = \omega$, then
for some $\alpha = \alpha(\langle f_\ell \such \ell < k \rangle)$ we 
have that:
\begin{equation}\label{flip}
\begin{split}
&\mbox{
For every  $u_\ell \in [\kappa \setminus \alpha]^{<\omega} $ and
$\eta_\ell \in {}^{u_\ell} 2$
the set}\\
& \{ n \in B \such (\forall \ell < k) (f^2_\ell(n) \cap 
\bar{A}^{[\eta_\ell]}
\neq \emptyset) \}\\
& \mbox{ is infinite},
\end{split}
\end{equation}
where  
\begin{equation*}
\begin{split}
\bar{A}^{[\eta_\ell]} & = \bigcap_{i \in u_\ell} A_i^{\eta_\ell(i)}, 
\mbox{ and }\\
A^\ell_i & = \left\{ \begin{array}{ll}
A_i, & \mbox{ if } \ell = 1,\\
\omega \setminus A_i, & \mbox{ if } \ell = 0.
\end{array}
\right.
\end{split}
\end{equation*}
\end{myrules}
\end{myrules}
\end{definition}

Remark: $f \in {\cal F}_g$ implies that $\bigcup_{i \in \dom(f)} f^2(i)$ 
is not $g$-slow.

\begin{claim}\label{3.3a}
We get an equivalent notion to ``$\bar{A}$ is $(g,\kappa)$-o.k.'',
if we modify the Definition~\ref{3.2}(c) as in (a) and/or 
as in (b), where
\begin{myrules}
\item[(a)]  We  demand \ref{3.2}(c) only for $f_\ell \in
{\mathcal F}_g$ that additionally satisfy
$\dom(f_0) = \cdots = \dom(f_{k-1}) = \omega$.
\item[(b)] We demand \ref{3.2}(c) only for $f_0, \dots , f_{k-1} 
\in {\mathcal F}_g$ such that
$\langle \min\{f^1_\ell(i) \such i < k \} \such i < B\rangle $ is
strictly increasing (we can even demand, increasing faster than
any given $h$), and 
for $i \in B$,
$\max\{f^1_\ell(i) \such \ell < k \} <
\min\{ f^1_\ell(i+1) \such \ell < k \}$.
\end{myrules}
\end{claim}

\proof
(a)
Suppose the $f_0,\dots , f_{k-1} \in {\mathcal F}_g$ in the original sense, 
and that we have required the analogue of \ref{3.2}(c) only for
${\mathcal F}_g$ in the restricted sense. We suppose that
$\bigcap_{\ell < k} \dom(f_\ell) = B$ and take a
strictly increasing enumeration $\{b_r \such r \in \omega\}$ 
of $B$.
Then we take $\tilde{f_\ell} \colon \omega \to 
[\omega]^{<\omega}$, $\tilde{f}_\ell(r) =
f_\ell(b_r)$ for $r \in 
\omega$.
The analogue of \ref{3.2} for the ${\mathcal F}_g$ 
in the restricted sense gives $\alpha \in \kappa$ and 
infinite intersections in \eqref{flip} for the $\tilde{f_\ell}$.
The intersections are
also infinite for the original $f_\ell$.

(b) Suppose that
 $k < \omega$, $f_0, \dots , f_{k-1} \in {\mathcal F}_g$,
$\bigcap_{\ell \in \omega} \dom(f_\ell) = B \in [\omega]^\omega$
and $\limsup\langle \min\{m_{f_\ell(i)} \such \ell \in k \}
\such i \in B \rangle = \omega$.
Then we can thin out the domain $B$ to some infinite $B'$,
inductively on $i$
such that the $f_\ell \restriction B'$ fulfil all the requirements
from \ref{3.3a}(b).

\begin{crucialfact}\label{crucial}
Let $g \in {}^\omega \omega$.
If $\gr < \kappa = \cf(\kappa)$ and 
if there is some $\bar{A}$ that is $(g,\kappa)$-o.k., then every
$\Delta_1^1$-real that is computable in every function
$g' \geq^\ast g$ is \wneeded for the refining relation.
\end{crucialfact}

\proof
Let $\mathcal R = \{ B_\alpha \such \alpha < |{\mathcal R}|  \}$
witness $\gr < \kappa$. 
The family $\bar{A}$ is refined by $\mathcal R$:
For $i<\kappa$ for some $\alpha_i < |{\mathcal R}| $ and $\nu(i) 
\in \{0,1\}$ we have that $B_{\alpha_i} \subseteq A_i^{\nu(i)}$.
Since $\kappa$ is regular and since $\gr < \kappa$, there
are  for some $\ell < 2 $ and some $\beta < |{\mathcal R}|$ such that
$$Y= \{ i < \kappa \such \nu(i) = \ell \wedge \alpha_i = \beta \}$$ 
is unbounded.
So $B_\beta \subseteq \bigcap_{i \in Y} A_i^{\nu(i)}$.
We claim that $B_\beta$ is not $g$-slow.
Why?
Otherwise we have $C= \{ n < \omega \such | B_\beta \cap g(n)| > n \}
\in[\omega]^\omega$, and we may take
$f \in {\mathcal F}_g$ such that $C=\dom(f)$,
$f^1(n)=n$ and $f^2(n)= B_\beta 
\cap g(n)$.
Take any $\alpha \in \kappa$. Then we take $u_0 $ such that $u_0
= \{\gamma\}$, $\gamma \in Y$, $\gamma > \alpha$ and 
$\eta_0 =\{(\gamma,0)\}$ 
and $\eta'_0 =\{(\gamma,1)\}$.
Then we do not have $(\exists^\infty n) f^2(n)
 \cap A_\gamma^0 \neq \emptyset$
and  
$(\exists^\infty n) f^2(n) \cap A_\gamma^1 \neq \emptyset$ at the same time,
because $B_\beta$ is refining $A_\gamma$. So $\bar{A}$ is not
$(g,\kappa)$-o.k., in contrast to our assumption.

But now we can compute recursively from $B_\beta$ some $g'
\geq^\ast g$, for example we may take $g'(n)=$(the $n$th 
element of $B_\beta$) $+1$.
Hence every hyperarithmetic real
that is computable in every function
$g' \geq^\ast g$  is recursive in $B_\beta$.
\proofend

So, how do we get the premises of the crucial fact?
The rest of this section will be devoted to this issue.
We take $g$ growing sufficiently fast so that
every $\Delta_1^1$-function is computable in every $g' \geq g$.
Such a $g$ exists by \cite{Jockusch, solovay:hyper} and the fact that there
are only countably many $\Delta^1_1$-functions.
We fix such a $g$.
We  consider the case $\kappa = \cf(\kappa) > \aleph_1$ and intend to 
show the consistency of ``$\gr = \aleph_1$ 
and there is some $\bar{A}$ that is $(g,\kappa)$-o.k.''

\begin{definition}\label{3.5}
\begin{myrules}
\item[(1)] $K_g=K= \{(P,\name{\bar{A}}) \such 
P$ is a ccc forcing and $ \Vdash_P$ 
``$\name{\bar{A}}$ is  $(g,\kappa)$-o.k.''$\}.$ 
For a fixed $g$, we often leave out the subscript.

\item[(2)] $(P_1, \name{\bar{A_1}}) \leq_K (P_2,\name{\bar{A_2}})$
iff $P_1 \lessdot P_2$  and $\name{\bar{A_1}} = \name{\bar{A_2}}$.
\end{myrules}
\end{definition}

\begin{claim}\label{3.6}
\begin{myrules}
\item[(1)]
We have that $K \neq \emptyset$. In fact,
if $P$ is the forcing adding $\kappa$ Cohen reals and $\name{\bar{A}}$
is the enumeration of the $\kappa$ Cohen reals,
then $(P,\name{\bar{A}}) \in Kg$ for any function $g$.
(This is true for any function $g$.)

\item[(2)] 
If $(P_\alpha, \name{\bar{A}}) \in K$ for 
$\alpha < \delta$, $\delta$ a limit cardinal, and
$\langle P_\alpha \such \alpha < \delta \rangle$
is increasing and continuous, and $P=\bigcup_{\alpha < \delta} 
P_\alpha$, then $(P,\name{\bar{A}}) \in K$ and $\alpha < \delta 
\Rightarrow (P_\alpha,\name{\bar{A}}) \leq _K
(P,\name{\bar{A}})$.

\nothing{\item[(3)]
If $(P,\name{\bar{A}})  \in K$ and $\Vdash_P
\mbox{``} \name{Q}$ is a forcing notion of cardinality $\mu < \kappa$'',
then $(P,\name{\bar{A}}) \leq_K (P \ast \name{Q},\name{\bar{A}}) \in K$.}
\end{myrules}
\end{claim}

\proof (1) Suppose that $f_0,\dots, f_{k-1}
\in V[G_\kappa]$ are injective functions. 
We take $\alpha$ such that $f_0,\dots, f_{k-1}
\in V[G_\alpha]$ where $G_\alpha$ is a generic filter for the 
first $\alpha$ Cohen reals. Suppose that $\eta_\ell \in {}^{u_\ell}
2$. Now a density argument gives that these 
 $\bar{A}^{[\eta_\ell]}$ 
``flip for infinitely many $n \in B$'' to 0 or to 1 within 
$f^2_\ell(n)$ for every $\ell < k$.

(2) $P$ has the c.c.c.\ by a Fodor argument.
Now we show that
$ \Vdash_P$ 
``$\name{\bar{A}}$ is  $(g,\kappa)$-o.k.''$\}.$ 
Only the case of $\cf(\delta) = \omega$ is not so easy.
We suppose that $\delta = \bigcup_{n \in \omega}  \alpha(n)$,
$0< \alpha(n)< \alpha(n+1)$.
Towards a contradiction we assume that $p^\ast \in P_{\alpha(0)}$, and
$$p^\ast \Vdash \mbox{``}\name{B},
\langle \name{f_\ell} \such \ell < k \rangle
\mbox{ form a counterexample to $\bar{A}$ being $(g,\kappa)$-o.k.''}
$$
For each $n \in \omega$ we find $\langle q_{n,i} \such i \in \omega
\rangle$ such that
\begin{myrules}
\item[$(\alpha)$] $q_{n,i} \in P$,
\item[$(\beta)$] $q_{n,0} = p^\ast$,
\item[($\gamma$)]
$P \models q_{n,i} \leq q_{n,i+1}$,
\item[$(\delta)$]  for some  $\name{b_{n,i}}$,
$\name{f^1_{n,\ell,i}}$, $\name{f^2_{n,\ell,i}}$
$P_{\alpha(n)}$-names we have 
$$
q_{n,i} \Vdash \mbox{``}\name{b_{n,i}}
\mbox{ is the $i$-th member of } \name{B},
\name{f_\ell}(\name{b_{n,i}}) =
(\name{f^1_{n,\ell,i}},\name{f^2_{n,\ell,i}})\mbox{''},
$$
\item[$(\eps)$] $q_{n,i} \restriction \alpha(n)
= q_{n,0} \restriction \alpha(n) =
p^\ast \restriction \alpha(n)$.
\end{myrules}

How do we choose these? 
Let $n$ and $\alpha(n)$ be given. Then we choose
$q'_{n,i}$ increasing in $i$ such that
$q'_{n,i} \in P$ and   $b'_{n,i}$,
$(f^1)'_{n,i}$, $(f^2)'_{n,\ell,i}$
in $V$ and
$$q'_{n,i} \Vdash \bigwedge_{\ell < k}
\mbox{ the $i$th element of $\name{B}$}
=\check{b'_{n,i}} \wedge \name{f_\ell}(\check{b'_{n,i}})=
(\check{(f^1)'_{n,\ell,i}},\check{(f^2)'_{n,\ell,i}}).$$
Then we take 
\begin{eqnarray*} 
\name{b_{n,i}} = (b'_{n,i}, q'_{n,i} \restriction P_{\alpha(n)}),\\
\name{f^1_{n,\ell,i}}=((f^1)'_{n,\ell,i}, q'_{n,i} 
\restriction P_{\alpha(n)}),\\
\name{f^2_{n,\ell,i}}=((f^2)'_{n,\ell,i}, q'_{n,i} 
\restriction P_{\alpha(n)}),\\
p_{n,i} = p^\ast\restriction \alpha(n) \cup 
q'_{n,i} \restriction [\alpha(n),\delta).
\end{eqnarray*}
Here, the restriction $\restriction \alpha$ is any reduction
function witnessing $P_\alpha \lessdot P$ (see \cite{abraham:handbook}),
and in the general case, if $P_\alpha$ is not
the initial segment of length $\alpha$ of some iteration,
the term $q'_{n,i} \restriction [\alpha(n),\delta)$ 
has to be interpreted as some element from a quotient forcing algebra.

\medskip

Now for every $n$ we define $P_{\alpha(n)}$-names
\begin{equation*}
\begin{split}
\name{B'_n} &= \{ \name{b_{n,i}} \such i < \omega \},\\
\name{f_{\ell,n}} &\colon \name{B'_n} \to V,\\
\name{f_{\ell,n}}(\name{b_{n,i}}) & =
(\name{f^1_{\ell,n}}(\name{b_{n,i}}),
\name{f^2_{\ell,n}}(\name{b_{n,i}}))=
(\name{f^1_{\ell,n,i}},
\name{f^2_{\ell,n,i}}).
\end{split}
\end{equation*}

Now we have that
\begin{equation*}
\begin{split}
p^\ast & \Vdash \mbox{``} \name{B'_n} \in [\omega]^{\aleph_0} ,
\name{f_{\ell,n}} \mbox{ is a function with domain } 
\name{B'_n} \mbox{ and }\\
&\mbox{$\limsup \langle \name{f^1_{\ell,n}}(b) \such b 
\in \name{B'_n} \rangle
=\omega$ and }\\ 
& \name{f^2_{\ell,n,i}} \mbox{ when defined is a subset of } [0,
g({\name{f^1_{\ell,n,i}}}))
\mbox{ of cardinality } > \name{f^1_{\ell,n,i}}\mbox{''}.
\end{split}
\end{equation*}

\smallskip

As $(P_{\alpha(n)}, \name{\bar{A}})$ is in $K$
we for every $n$
\begin{equation*}
\begin{split}
p^\ast\restriction \alpha(n)  & \Vdash_{P_{\alpha(n)}} \mbox{`` for some 
$\name{\beta} < \kappa$ for every $u_\ell \subseteq
[\kappa \setminus \name{\beta}]^{\aleph_0}$ for
every $\eta_\ell \in {}^{u_\ell} 2$}\\
& \left\{ b \in \name{B'_n} \such \bigwedge_{\ell < k }
\name{{f^2_{\ell,n}(b)}} \cap \name{\bar{A}}^{[\eta_\ell]} 
\neq \emptyset \right\} \mbox{ is infinite.''}
\end{split}
\end{equation*}

Let $\name{\beta_n} < \kappa$ be such a $P_{\alpha(n)}$-name.
Since $P_{\alpha(n)}$ has the ccc, 
there is some $\beta^\ast_n < \kappa$ such that
$\Vdash_{P_{\alpha(n)}} \name{\beta_n} < \beta^\ast_n < \kappa$.
Since $\kappa$ is regular we have that $\beta^\ast = 
\bigcup_{n \in \omega} \beta^\ast_n < \kappa$.

It suffices to prove that
$$p^\ast \Vdash \mbox{``$\beta^\ast$ is as required in the definition of
$(g,\kappa)$-o.k.''}
$$

If not, then there are counterexamples
$u_\ell \in [\kappa \setminus \beta^\ast]^{<\aleph_0}$,  $
\eta_\ell \in {}^{u_\ell}2$, $q$ and $b^\ast$ such that

\begin{equation}\tag{$\diamond$}\label{contr}
\begin{split}
p^\ast & \leq q \in P = P_\delta\\
q & \Vdash \mbox{``}\left\{b \in \name{B} \such
(\forall \ell < k) (\name{f^2_\ell}(b) \cap \name{\bar{A}^{[\eta_\ell]}}
\neq \emptyset) \right\} \subseteq [0,b^\ast]\mbox{''.}
\end{split}
\end{equation}

For some $n(\ast)< \omega$ we have that $q \in P_{\alpha(n(\ast))}$.
Let $G \subseteq P$
be generic over $V$, and let $q \in G_{\alpha(n(\ast))}$.
So by the choice of ${\beta_{n(\ast)}} <
\beta^\ast$ we have that
$$
p \Vdash_{P_{\alpha(n(\ast))}} C = \{ b \in \name{B'_{n(\ast)}} \such
(\forall \ell < k)
(\name{f^2_{\ell,n(\ast)}}(b)
\cap \name{\bar{A}^{[\eta_\ell]}} \neq \emptyset) \}
\mbox{ is infinite''.}
$$
Recall that $\name{B'_{n(\ast)}}$ and
$\name{f_{\ell,n(\ast)}}(b)$ are 
$P_{\alpha(n(\ast))}$-names and that 
$\name{\bar{A}^{[\eta_\ell]}}$ is a $P_0$-name.

Now $\name{B'_{n(\ast)}}= \{ \name{b_{n(\ast),i}} \such i < \omega \}$,
so for some $i$ we have that $\name{b_{n(\ast),i}}[G] > b^\ast$.
So $q_{n(\ast), i}\in G\cap P_{\alpha(n(\ast))}$ forces
``the $i$-th member of $\name{B}$ is $\name{b_{n(\ast),i}}$
and $\name{f_\ell}(\name{b_{n(\ast),i}})=
\name{f_{\ell,n(\ast)}} (\name{b_{n(\ast),i}}) =
(\name{f^1_{\ell,n(\ast),i}}, \name{f^2_{\ell,n(\ast),i}})$.
Note that $q_{n(\ast),i} \restriction \alpha(n(\ast))
= p^\ast \restriction \alpha(n(\ast))$ according to $\eps)$, and
hence $q_{n(\ast),i} \not\perp q$.
So 
there is some $r \geq
q$ and $r \geq  q_{n(\ast),i}$. Such an $r$
forces the contrary of the property forced in \eqref{contr}, and 
finally we reached a contradiction.
\nothing{
(3) Take a bijection of $\mu \times \kappa$ to $\kappa$ and
treat the names for the function $f_\ell$ by allowing $\mu$ places 
for each value.}
\proofend

Now \ref{3.7} and \ref{3.8} are like \cite{Sh:707}. 
For $h \colon \omega \to \omega$ We write 
$\lim_D \langle h(i)\such i \in \omega \rangle =
\omega$ if for all $m <\omega$ we have that
$\{i \such h(i) > m\} \in D$.

\begin{claim}\label{3.7}
Assume that in $V$:
\begin{myrules}
\item[(a)] $\bar{A}$ is $(g,\kappa)$-o.k.
\item[(b)] $\kappa = 2^{\aleph_0}$.
\end{myrules}
Then there is an ultrafilter $D$ on $\omega$ such that
\begin{equation}\tag{$\ast$}\label{goodultrafilter}
\begin{split}
& \mbox{if $f \in {\mathcal F}_g$ and $\dom(f) \in D$
and $\lim_D\langle f^1(i) \such i \in \dom(f) \rangle = \omega$}\\
& \mbox{then for some $\alpha_f
< \kappa$ for every $u \in [\kappa \setminus \alpha_f]^{<\aleph_0}$
 and $\eta \in {}^u 2$}\\
& \mbox{we have that $\{ n \in \dom(f) \such f^2(n) \cap \bar{A}^{[\eta]}
\neq \emptyset \} \in D$.}
\end{split}
\end{equation}
\end{claim}

\nc{\AP}{{\mathcal{AP}}}

\proof
Let ${\mathcal F}_g = \{ f_j \such j < \kappa \}$.
Let $\AP$ be the set of tuples $(D,i,\alpha)$ such that

\begin{myrules}
\item[(i)] $D$ is a filter on $\omega$ containing the co-finite
subsets, $\emptyset \not\in D$, $i,\alpha < \kappa$,
\item[(ii)] $D$ is generated by $< \kappa$ members,
\item[(iii)] if $k < \omega$ and for $\ell <k$,
$j_\ell < i$, and $\dom(f_{j_\ell}) \in D$
and $\lim_D\langle f^1_{j_\ell}(i) \such i \in \dom(f_{j_\ell})
\rangle = \omega$ and $u_\ell
\in [\kappa \setminus \alpha]^{<\aleph_0}$, 
$\eta_\ell \in {}^{u_\ell} 2$, then
$$ \left\{ n \in \bigcap_{\ell < k} \dom(f_{j_\ell})
\such \bigwedge_{\ell < k} \Bigl(f^2_{j_\ell}(n) 
\cap \bar{A}^{[\eta]} \neq
\emptyset \Bigr) \right\}
\neq \emptyset {\text{ mod }} D.$$
\end{myrules}
Let $(D_1,i_1,\alpha_1) \leq_{\AP} (D_2, i_2,\alpha_2)$ if both
tuples are in $\AP$ and 

\begin{myrules}
\item[($\alpha$)] $D_1 \subseteq D_2$, $i_1 \leq i_2$, 
$\alpha_1 \leq \alpha_2$, and
\item[$(\beta)$] if $k < \omega$ and $\{j_0, \dots, j_{k-1}\}
\subseteq i_1$, $\dom(f_{j_\ell}) \in D_2$
and $\lim_{D_2}\langle f^1_{j_\ell}(i) \such i \in \dom(f_{j_\ell})
\rangle = \omega$ and $u_\ell \subseteq [\alpha_1,\alpha_2)$ is finite and 
$\eta_\ell \in {}^{u_\ell} 2$ then
$$\left\{ n \in \bigcap_{\ell < k} \dom(f_{j_\ell}) \such
\bigwedge_{\ell < k}
f^2_{j_\ell(n)} \cap \bar{A}^{[\eta_\ell]} \neq \emptyset \right\} 
\in D_2.$$
\end{myrules}

Now we have that
\begin{myrules}
\item[$\boxtimes_1$] $(\AP,\leq_{\AP}) $ is a non-empty partial order.
Take $i=\alpha=0$ and $D$ the filter of all cofinite subsets of $\omega$.

\item[$\boxtimes_2$] In $(\AP,\leq_{\AP}) $ every increasing sequence
of length $< \kappa$ has an upper bound, namely, take 
the filter generated by the union in the first coordinate and
take the supremum in the second and in the third coordinate.
\item[$\boxtimes_3$]
If $B \subseteq \omega$ and $(D,i,\alpha) \in \AP$ then there are some
$D'$, $i'$, $\alpha'$ such that
$(D',i',\alpha') \geq_{\AP} (D,i,\alpha)$ 
and that $B \in D'$ or that $\omega \setminus B \in D'$.
Why? Try $D' =$ the filter generated by 
$D \cup \{B\}$ and the same $i$ and $\alpha$.
 If this fails then we can find $k < \omega$,
such that for $\ell < k$ we have  $j_\ell < i$,
such that $\dom(f_{j_\ell}) \in D'$
and $\lim_{D'}\langle f^1_{j_\ell}(i) \such i \in \dom(f_{j_\ell})
\rangle = \omega$,
$u_\ell \in [\kappa \setminus \alpha]^{<\aleph_0} $, $\eta_\ell
\in {}^{u_\ell} 2$ and such that
$$\left\{ n \in 
\bigcap_{\ell < k} \dom(f_{j_\ell}) \such 
f^2_{j_\ell}(n) \cap \bar{A}^{[\eta_\ell]} \neq 
\emptyset \right\} \cap B = \emptyset {\text { mod }} D.
$$
Let $\alpha' < \kappa$ be such that $\alpha \leq \alpha'$ and 
$\bigwedge_{\ell < k} u_\ell \subseteq\alpha'$.
Let $D'$ be the filter generated by 
\begin{multline*}
\makebox[1cm]{} D \cup \biggl\{ \Bigl\{ n \in \bigcap_{\ell < k } 
\dom(f_{j_\ell}) \such
f^2_{j_\ell}(n) \cap \bar{A}^{[\eta_\ell]} \neq\emptyset \Bigr\}
\such \\
k < \omega, j_\ell < i, u_\ell \in 
[\alpha'\setminus \alpha]^{<\aleph_0},
\eta_\ell \in {}^{u_\ell} 2 \biggr\}.
\end{multline*}
Then $\omega \setminus B \in D'$, and $(D',i,\alpha') \in \AP$.

\item[$\boxtimes_4$] 
If $(D,i,\alpha) \in \AP$  then for some $D'$, $\alpha'$ we have that
$(D', i+1, \alpha') \in \AP$.
\proof
Let $M \prec (H(\chi),\in)$ such that $M \cap \kappa 
\in \kappa$, $(D,i,\alpha) \in M$, ${\mathcal F}_g \in M$, and 
$|M| < \kappa$. 
Suppose that $\dom(f_i) \in D$ and that
$\lim_D \langle f^1_i(k) \such k \in \dom(f_i)\rangle = \omega$.
Let $\alpha' = M \cap \kappa$.
Let $D_1$ be the filter in the boolean algebra in 
${\mathcal P}(\omega) \cap M$ generated by
\begin{multline*}
\makebox[1cm]{} (D\cap M) \cup \biggl\{ 
\Bigl\{ n \in \bigcap_{\ell < k } \dom(f_{j_\ell}) \such
f^2_{j_\ell}(n) \cap \bar{A}^{[\eta_\ell]} \neq\emptyset
\Bigr\}
\such \\
k < \omega, j_\ell \leq i, u_\ell \in 
[\alpha'\setminus \alpha]^{<\aleph_0},
\eta_\ell \in {}^{u_\ell} 2 \biggr\}.
\end{multline*}
Since in $M$,
$\bar{A}$ is $(g,\kappa)$-o.k., this has the infinite intersection 
property.
Let $D_2'$ be an ultrafilter in $M$ extending $D_1$.
Let $D'$ be the filter on $\omega$ in $V$ that $D_2'$ generates.
\end{myrules}

Now we take a maximal element in the partial order
$(\AP,\leq_{\AP}) $. By the properties
$\boxtimes_1$ to $\boxtimes_4$ it is as required in 
\eqref{goodultrafilter}.
\proofend

Note that $(\ast)$ of \ref{3.7} implies that $\bar{A}$ is $(g,\kappa)$-o.k.
The following is a preservation theorem
for suitable ultrafilters:

\begin{claim}\label{3.8}
Assume that 
\begin{myrules}
\item[(a)] $\bar{A}$ is $(g,\kappa)$-o.k.
\item[(b)] $D = \langle D_\eta 
\such \eta \in {}^{<\omega}\omega \rangle$, $D_\eta = D$,
$D$ is  ultrafilter on $\omega$ as in \ref{3.7}.
\item[(c)] $ Q_D=\{ T \such T \subseteq {}^{<\omega}\omega $ 
is a subtree, and for some $\eta \in T, 
\eta \trianglelefteq \nu \in T 
\Rightarrow \{k \such \nu \concat k \in T \} \in D_\nu \}$, 
ordered by inverse inclusion. (The $\triangleleft$-minimal $\eta$ of this
sort is called the trunk of $T$, $\trunk(T)$.)
\end{myrules}
Then $\Vdash_{Q_D} \mbox{``} \bar{A} \mbox{ is } (g,\kappa)\mbox{-o.k.''}$.
\end{claim}

\proof We use the fact \cite{Sh:707} that
$Q_D$ has the pure decision property:
Let $\varphi_i$, $i \in \omega$, be countably many sentences of
the $Q_D$-forcing language. We think of names $\name{f_\ell}$,
$\ell < k$, for some
elements of ${\mathcal F}_g$ and $\varphi_i =$
``$\Bigl($the $i$-th element of $\name{B}
=\bigcap_{\ell < k}\dom(\name{f_\ell})\Bigr)
= \check{b_i}$
and $\bigwedge_{\ell<k}
\name{f_\ell}(\check{b_i}) = (\check{f^1_{\ell,i}},\check{f^2_{\ell,i}})$''.
The pure decision property says:

\begin{equation*}
\forall p \in Q_D
\; \exists q \geq_{tr} p \;
\forall r \geq q \; \forall i\; 
\Bigl(r \Vdash \varphi_i \rightarrow (\exists s_i \in r) 
 q^{[s_i]} \Vdash \varphi_i \Bigr),
\end{equation*}

where we write $\geq_{tr}$ for the pure extension:
$q \leq_{tr} r$ if $r \subseteq q$ and $\trunk(q)=\trunk(r)$,
and $q^{[s_i]} = \{ \eta \in q \such s_i \trianglelefteq \eta \}$.

Towards a contradiction we 
assume that there is a counterexample.
By Claim~\ref{3.3a} (first (b) and then (a))
we may assume that
it is of the following form
\begin{equation*}\label{astast}\tag{$\ast\ast$}
\begin{split}
p^\ast \Vdash & \mbox{``}
\langle \name{f_\ell} \such \ell < k \rangle
\mbox{ form a task} \\
& \mbox{ such that the intersection of the
domains is $B = \omega$}\\
& \mbox{ and for $i \in B$, $\max\{\name{f^1_\ell}(i) \such \ell < k \} 
<\min\{ \name{f^1_\ell}(i+1) \such \ell < k \}$}\\
&\mbox{ and there is no $\alpha < \kappa$
such that the statement}\\
&\mbox{  \eqref{flip} from Definition \ref{3.2}(3)(c) holds.''}
\end{split}
\end{equation*}

\nothing{ What if $\bigcup_{\ell<k} \dom(\name{f_\ell}) \not\in D$}
\nothing{
So we have for some $m^\ast \in \omega$ that
\begin{equation*}
\begin{split}
p^\ast \Vdash & \mbox{``}\name{B}=\bigcap_{\ell \in k} \dom(\name{f_\ell})
\mbox {is infinite, and } \forall \alpha \in \kappa \;
\exists u_\ell \in [\kappa \setminus \alpha]^{<\aleph_0} \;\\
&
\exists \eta_\ell \in {}^{u_\ell} 2
\{ n \in \name{B} \such (\forall \ell < k ) \name{f_\ell}(n) \cap
\bar{A}^{[\eta_\ell]} \neq \emptyset \} \subseteq [0,m^\ast].
\end{split}
\end{equation*}
}

We find $q$ such that
\begin{myrules}
\item[$(\alpha)$] $q \in P$
\item[$(\beta)$] $q \geq_{tr} p^\ast$,\\
\item[($\gamma$)]
 for  all $i \in \omega$ for all $f^1_{\ell,i} \in \omega$, $f^2_{\ell,i}
\subseteq [0,g(f^1_{\ell,i}))$ of size bigger than $f^1_{\ell,i}$
we have that
\begin{equation*}
\begin{split}
&\mbox{ if } r \geq q, r \Vdash \mbox{``}
\name{f_\ell}(\check{i}) =
(\check{f^1_{\ell,i}}, \check{f^2_{\ell,i}})\mbox{''},\\
&\mbox{then also
for some $s_i \in r$, the condition
 $q^{[{s_i}]}$ forces the same.''}
\end{split}
\end{equation*}
\end{myrules}
We fix such a $q$.

Now we set for $\nu \in q$ and $\ell < k$
$$B^1_{\nu,\ell} = \{ i \in \omega \such
\mbox{ some pure extension of } q^{[\nu]} \mbox{ decides } \name{f_\ell(i)}
\}.$$

\smallskip

We say $(\nu,\ell)$ is 1-good if $B^1_{\nu,\ell} \in D$. Let
for $i \in B^1_{\nu,\ell}$,
$h_{\nu,\ell}(i)=(h^1_{\nu\ell}, h^2_{\nu,\ell}) $ 
the value of $\name{f_{\ell}(i)} $ that is given 
by the pure decision. This is well-defined because
any two pure extensions are compatible.
Of course, 
by the requirements we had put  on the
counterexample,  we have that $\lim_D\langle
h^1_{\nu,\ell}(i) \such i \in B^1_{\nu,\ell} \rangle = \omega$.

We say that $(\nu,\ell) \in q \times k$ is 2-good, if it is not 1-good and
we have for all $m \in \omega$ that 
\begin{equation*}
\begin{split}
M_{\nu,\ell,m} =
\left\{ 
j \in \omega \such \right. & 
(\exists i \in \omega) 
(h_{\nu\concat j,\ell}(i))
\mbox{ is well-defined,}\\
&
 \left.\mbox{ and } h^1_{\nu\concat j,\ell}(i)> m) \right\} \in D.
\end{split}
\end{equation*}

So, for 2-good but not 1-good 
$(\nu,\ell)$ we may define for $j \in M_{\nu,\ell,m}$, 
\begin{equation*}
\begin{split}
g_{\nu,\ell}(j)= &
h_{\nu\concat j,\ell}(i_{\nu\concat j,\ell}),\\
&
\mbox{ where $i_{\nu\concat j,\ell}$ is 
such that $h_{\nu\concat j,\ell}(i_{\nu\concat j,\ell})$
is defined in $h^1_{\nu\concat j,\ell}(i_{\nu\concat j,\ell})>m$}\\
& \mbox{ and if there is a maximal such $i$, then take this as
 $i_{\nu\concat j,\ell}$.}
\end{split}
\end{equation*}

We show that there is $M'_{\nu,\ell,m}\in D$, $M_{\nu,\ell,m}
\supseteq M'_{\nu,\ell,m}$ such that for $j \in M'_{\nu,\ell,m}$
there a maximal such $i$: 
If $h_{\nu\concat j,\ell}(i)$ is defined and $i' < i$ then there is some pure
extension deciding $h_{\nu\concat j,\ell}(i')$ since 
there are only finitely many
possibilities for it values, by the third line of \eqref{astast}.
Hence some pure extension decides the value.
Hence also $h_{\nu\concat j,\ell}(i')$ is defined. If 
$h_{\nu\concat j,\ell}(i)$ is
defined for all $i$, then $(\nu\concat j,\ell) $ is 1-good.
Hence , if $(\nu\concat j, \ell)$ is 2-good but not 1-good,
then there is a maximal $i$ witnessing
$j \in M_{\nu,\ell,m}$.
If $\{j \such (\nu\concat j,\ell) $ is 1-good $\} \in D$, then
by gluing together suitable pure extensions $r_j$ of
 $q^{[\nu\concat j]}$ together we get a pure extension of
$q^{[\nu]} $  that shows that 
$(\nu,\ell)$ is 1-good. Hence $X= \{ j \such 
(\nu\concat j $ is 2-good and not 1-good $\} \in D$. 
So we may take  $M'_{\nu,\ell,m}
= M_{\nu,\ell,m}\cap X$.
In order to simpily notation, we  assume that  $M'_{\nu,\ell,m}
= M_{\nu,\ell,m}$.

\smallskip

Also from the third line of \eqref{astast} we get
that for every $\nu \in q$ either for 
all $\ell<k$, $(\nu,\ell) $ is 1-good or
no $(\nu,\ell)$ is 1-good. In the latter case there 
is some $i_\nu$, such that
for all $\ell<k$, 
$\dom(h_{\nu,\ell}) = i_\nu$ or $\dom(h_{\nu,\ell})=i_\nu +1$.
Moreover, also by \eqref{astast} we get that
if for some $\ell < k$, for all $m$, $M_{\nu\ell,m} \in D$, then
for all $\ell <k$, for all $m$, $M_{\nu,\ell,m} \in D$. 
So if $(\nu,\ell)$ is 2-good, then all $(\nu,\ell')$ are 2-good.
We call $\nu$ $i$-good if there is some $\ell$ such that $(\nu,\ell)$ is
$i$-good. We set $M_{\nu,m} = \bigcap_{\ell < k} M_{\nu,\ell,m}$. 

\smallskip

We fix some diagonal intersection $M_{\nu}$ of
$\langle M_{\nu,m} \such m \in \omega \rangle$, such that
$\lim\langle i_{\nu\concat j} \such 
j \in M_{\nu} \rangle = \omega$.

Then we also have that
$\lim_D\langle \min\{g^1_{\nu,\ell}(j) \such \ell < k \} \such 
j \in M_{\nu}\rangle = \omega$, because for each $z< \omega$,
$\{ j \such  \min\{g^1_{\nu,\ell}(j) \such \ell < k \} < z \}$
is a cofinite set. Hence
 $g_{\nu,\ell} \in {\mathcal F}_g$.
By combining with an enumeration of $M_\nu$, we may assume
that $\dom(g_{\nu,\ell}) = \omega \in D$.
We will not write this enumeration, in order to prevent too
clumsy notation, but
we shall later apply that $D$ is as in \ref{3.7}
for ${\cal F}_g$, and therefore we need that the domains are in $D$.

\smallskip

Now we take $\chi$ sufficiently large and $N \prec (H(\chi),\in)$ such that
$\langle \name{f}_\ell \such \ell < k \rangle \in N$,
$\langle B^1_{\nu,\ell}, h_{\nu,\ell},
g_{\nu,\ell} \such \nu \in q, \ell < k \rangle \in N$, $q, D  
\in N$.  We take $\alpha^\ast = \sup(N \cap \kappa)$.
We claim that $q$ forces that $\alpha^\ast$ is as in the 
Definition~\ref{3.2}(3)(c).

\smallskip

\nothing{
Now we define for every branch
$\bar{s} = \langle s_i \such i \in \omega \rangle $ of $q$
(the $s_i$ need not be of length $i$, we just require $s_i 
\triangleleft s_{i+1}$) the set $B^{\bar{s}} \in V$ by
\begin{equation*}
\begin{split}
B^{\bar{s}} &= \{m_i \such m_i \mbox{ as above and this is forced
by } q^{[s_i]} \},\\
f^{\bar{s}}_\ell&\colon B^{\bar{s}} \to V,\\
f^{\bar{s}}_\ell(m_{i}) & = w_{\ell,i} \mbox{ if $s_i$ forces this }.
\end{split}
\end{equation*}
Now we take $p^{\ast\ast}$ such that $p^{\ast\ast} \geq q_n$ for all
$n$ and such that there is some infinite set $B''\subseteq B'$ such that
for $n \in B''$ we have that $s_n \in p^{\ast\ast}$.
This ist possible by taking appropriate pseudointersections. $D$ is
a $p$-point.
Now we take only countably many branches $\{\bar{s}^k \such k \in \omega \}$
which are dense in $q$ in the following sense:
$(\forall s \in q)(\exists k \in \omega)(\exists i) s \triangleleft s^k_i$.
Now we have, because $\bar{A}$ is $(g,\kappa)$-o.k.\ in $V$ and $D$ is
as in \ref{3.6}, that
\begin{equation*}\tag{$\diamond$}
\label{contr1}
\begin{split}
& \mbox{for all $k$,  for some 
$\beta_k < \kappa$ for every $u_\ell \subseteq
[\kappa \setminus \beta_k]^{\aleph_0}$ for
every $\eta_\ell \in {}^{u_\ell} 2$ we have that}\\
& C^{\bar{s}^k}_{\eta_\ell} = \left\{ m \in B^{\bar{s]^k}} 
\such \bigwedge_{\ell < k }
f^{\bar{s}^k}_\ell(m) \cap \bar{A}^{[\eta_\ell]} 
\neq \emptyset \right\}\in D \mbox{.}
\end{split}
\end{equation*}
We take $\beta= \bigcup_{k \in \omega } \beta_k$.
It suffices to prove that
$$q \Vdash \mbox{`` $\beta$ is as required 
in the definition of $(g,\kappa)$-o.k.''}
$$
Since we chose a counterexample according to Claim~\ref{3.3a},
we have for every $\nu \in q$:
All $f^1_{\nu,\ell}$, $\ell < k$ have domain $\omega$ or 
for some $i_\nu \in \omega$, all $f^1_{\nu,\ell}$, $\ell < k$,
have domain $[0,i_\nu)$ or $[0,i_\nu]$.
}

If not, then there are counterexamples
$u_\ell \in [\kappa \setminus \alpha^\ast]^{<\aleph_0}$ and $
\eta_\ell \in {}^{u_\ell}2$ and $r \in Q_D$, $r \geq q$,
and $b^\ast$ such that

\begin{equation}\tag{$\diamond\diamond$}\label{contr2}
\begin{split}
r & \geq q, \mbox{ and } \\
r & \Vdash_{Q_D} \mbox{``}
\bigcap_{\ell < k} \dom(\name{f_\ell}) = \omega \mbox{ and }\\
&
(\forall i \in \omega)\max\{ \name{f^1_\ell}(i) \such \ell < k \}
< \min\{ \name{f^1_\ell}(i+1) \such \ell < k \}\\
& \mbox{ and }\left\{b \in \omega \such
(\forall \ell < k) (\name{f^2_\ell}(b) \cap \name{\bar{A}^{[\eta_\ell]}}
\neq \emptyset) \right\} \subseteq [0,b^\ast]\mbox{''.}
\end{split}
\end{equation}

\nothing{
Hence we may define for $\nu \in r$:
$\nu$ is 1-good if $i_\nu=\omega$, $\nu$ is 2-good if
$\omega = \lim_D \langle 
i_{\nu\concat j} \such j \in M\rangle$, but $\nu$ is not 1-good.
}

\smallskip

First case: 
There is some $\nu \in r$ with $\trunk(r) \trianglelefteq \nu$ 
such that all $\nu$ is 1-good. 
Now we take for each $t\in \omega$, 
some pure extension of $q_t^{[\nu]} $ of $r^{[\nu]}$ such that 
it forces 
$\bigwedge_{\ell <k} (h_{\nu,\ell} \restriction t
= \name{f_\ell}\restriction t)$.
Since $\bar{A}$ is $(g,\kappa)$-o.k., and since all is reflected to 
$N$ and by the choice of $\alpha^\ast$ we have that
$I=\{ n \in \omega \such (\forall \ell <k) (h^2_{\nu,\ell}(n) \cap
\bar{A}^{[\eta_\ell]} \neq \emptyset\}$ is infinite.
 So we take $t\in I$ such that $t> b^\ast$. Now $q_t^{[\nu]}$ contradicts 
\eqref{contr2}.

\smallskip

Second case. There is some  $\nu\in r$ such that all $\nu$,
$\ell < k$ are 2-good but not 1-good.
\nothing{
For 2-good  $\nu$ we choose $m(\nu,j,\ast)$ such that
\begin{myrules}
\item[$(\alpha)$] $m(\nu,j,\ast) < i_{\nu\concat j}$,
\item[$(\beta)$] under $(\alpha)$, $m(\nu,j,\ast)>j$ is maximal 
\end{myrules}}
We set $g_{\nu,\ell}(j) = h_{\nu\concat j,\ell}(i_{\nu\concat j,\ell})$ 
as purely
decided above $q^{[\nu\concat j]}$.
Fact: Now $\langle g_{\nu,\ell} \such \ell < k \rangle$ is as required in the 
definition of $\bar{A}$ being $(g,\kappa)$-o.k.,
because $\omega = \lim_D \langle g^1(i_{\nu\concat j }) 
\such j \in \omega \rangle$.
\nothing{
(For every $s < \omega$ we have for $D$-many $j$ that $i_{\nu\concat j} > s$
hence $m(\nu,j,\ast) > s$ but then $m_{g^\nu_\ell(j)} > s$.
}

Now we take for each $t\in \omega$, 
some pure extension of $q_t^{[\nu\concat j]} $ of $r^{[\nu\concat j]}$ 
such that 
it determines
$\bigwedge_{\ell <k} g_{\nu,\ell} \restriction t$.
Since $\bar{A}$ is $(g,\kappa)$-o.k., and since all is reflected to 
$N$ and by the choice of $\alpha^\ast$ we have that
$J=\{ n \in \omega \such (\forall \ell <k) (g^2_{\nu,\ell}(n) \cap
\bar{A}^{[\eta_\ell]} \neq \emptyset\}$ is infinite.
Then also $\hat{J}
= \{i_{\nu\concat n} \such n \in J\}$ is infinite.
 So we take $t> b^\ast$, $t \in \hat{J}$. Now the gluing together of 
$q_t^{[\nu\concat j]}$, $j \in \bigcap_{\ell < k}M_{\nu,\ell,t}$, 
contradicts 
\eqref{contr2} because we have $g_{\nu,\ell}(j)=
h_{\nu\concat j,\ell}(i_{\nu\concat j,\ell}) = 
f_\ell(i_{\nu\concat j})$, if $q_t^{[\nu]} \in G$.
Here we write $f_\ell$ for $\name{f_\ell}[G]$.

\nothing{
So assume that $r \geq q$ and $r \Vdash $ there is a counterexample 
with $\alpha^\ast$: $u_\ell \in [\kappa\setminus \alpha^\ast]^{<\aleph_0}$,
$\eta_\ell \in {}^{u_\ell} 2$, $\ell < k$ 
and $b^\ast$ form a counterexample (with the same earlier mentioned
$\name{f_\ell}$'s.
}

\smallskip

Third case: All $\nu \in r$ are neither 1-good nor 2-good. 
We shall prove something stronger:

\begin{quote}
An end-segment of the generic
$\bigcup \{\eta \such $ there is some element $q \in G$ with trunk $\eta\}$ 
can be thinned out (such that still infinitely 
many points are left) and injected into an infinite subset of
$\{ n \in \omega \such \bigwedge_{\ell < k} \name{f^2_\ell}[G](n)
\cap A^{[\eta_\ell]} \neq \emptyset \}$.
\end{quote}
This is more than enough. 

Let 
$i_{\nu,\ell}
= \max(B^1_{\nu,\ell}) < \omega$, because
$\nu$ is not 1-good.
Let $i^\ast_\nu = \dom(h_{\nu,\ell})$ such that
$i^\ast_\nu = i^\ast_{\nu,\ell}$ or $i^\ast_\nu =
i^\ast_{\nu,\ell} +1$. By the premise \eqref{astast},
there are such  $i^\ast_{\nu}$.
There is $r \geq q$ with no $\nu \in r$ being 1-good or 2-good in 
$N$. W.l.o.g. we take $q$ like that. Now we try to shrink $q$ purely.
Let $\nu_0 = \trunk(q)$. 

First: We have that
$\name{f_\ell}\restriction i^\ast_\nu$ is decided by $q$.
The range of 
$\langle i^\ast_{\nu \concat j} \such \nu \concat j \in q \rangle$ 
is bounded modulo $D$ because $\nu$ is not 2-good. Hence we
may assume that there is just one value $i_\nu^{\ast\ast}$.
So say (after shrinking $q$) that it is constant with value 
$i^{\ast\ast}_\nu \geq i^\ast_\nu$. 

Second we have that $\nu_0 \trianglelefteq \nu \in q$ 
implies that $q^{[\nu]}$ decides $\name{f_\ell}\restriction i^{\ast\ast}_\nu$.

Third we have that if $i \in [i^\ast_\nu, i^{\ast\ast}_\nu]$ then $\lim_D 
\langle f^1_{\nu\concat j,\ell}(i) \such j \in \omega \rangle= \omega$ by
the definition of $i_\nu^\ast$ and $i^{\ast\ast}_\nu$.
So define $g_{\nu,\ell,i}$ by $g_{\nu,\ell,i}(j) =
h_{\nu\concat j,\ell} (i)$.
So $g_{\nu,\ell,i} \in N$ is a function of the right form.
%

We have 
by the definition of $\alpha^\ast$ that for all $i \in 
[i^\ast_\nu,i^{\ast\ast}_\nu)$ for all $\nu \in q$ for all
$u_\ell$, $\eta_\ell$ that

$$A:=\{ b \such (\forall \ell < k)
g^2_{\nu,\ell,i}(b) \cap \bar{A}^{[\eta_\ell]} \neq \emptyset \}
\in D.$$
\nothing{
Now we take $j$ big enough to find some $\nu\concat j \in r$, $b
\in H_j$ such that $b > b^\ast$, and thus get a contradiction to
Equation~\eqref{contr2}.}
Since the range of $\bigcup\{\eta \such $ there is 
some element $q \in G$ with trunk $\eta\}=: \eta_\omega$
 is eventually contained is every set in $D$, we now
find the following infinite set:
We take $\langle \eta_n\such n \in \omega \rangle$ such that
$\eta_n \in \rge(\eta_\omega)\cap A$ and such that
$i^{\ast\ast}_{\eta_n} < i^\ast_{\eta_{n
+1}}$.
We set $\xi_n = \eta_n \restriction |\eta_n -1|$.
Then we have for almost all $n$ such that $\xi_n \in A$ and hence
for all $i \in [i^\ast_{\xi_n}, 
i^{\ast\ast}_{\xi_n})$:
$g_{\xi_n,\ell,i}(\eta_n(|\eta_n -1|)) =
h_{\xi_n \concat \eta_n(|\eta_n -1|),\ell}(i)=
h_{\eta_n,\ell}(i) = f_\ell(i)$.
So $\bigcup_{n \in \omega} [i^\ast_{\xi_n},i^{\ast\ast}_{\xi_n}]
\subseteq^\ast  \{ b \such (\forall \ell < k)
f^2_{\ell}(b) \cap \bar{A}^{[\eta_\ell]} \neq \emptyset \}
$ is infinite.

\proofend

\nothing{
$$\name{\eta}(t) \in 
\{ n \in \omega \such \bigwedge_{\ell < k} w_{\name{f_\ell}(n)}
\cap A^{[\eta_ell]} \neq \emptyset \}.$$
But now, according to equation
\eqref{contr1}, there is $m_i > m^\ast$, such that
$s_i^k \in r$ amd $m_i 
\in C^{\bar{s}}_{\eta_{\ell}}$.
We have that
 \begin{equation*}
\begin{split}
q^{[s_i]} &\Vdash 
\check{m_i} \in \name{B}  
\;\wedge \; \forall \ell < k \name{f}_\ell(\check{m_i})
\cap \bar{A}^{[\eta_\ell]} \neq \emptyset, \mbox{ and}
\\
r^{[s_i]} & \geq r \mbox{ and } r^{[s_i]} \geq q^{[s_i]},
\end{split}
\end{equation*}
so a contradiction to \eqref{contr2}.
}

\begin{claim}\label{3.9} Let $\kappa = \cf(\kappa) > \omega_1$.
Let ${\bf V}_0 \models \CH$ and let $P_0 = {\mathbb C}_\kappa$ be the forcing 
adding $\kappa$ Cohen reals.
We fix some function $g \in {\bf V}_0$, so that
every hyperarithmetic function in ${\bf V}_0$
is computable in every $g' \geq g$.
Set ${\bf V}_1 = {\bf V}_0[G_0]$. Let in ${\bf V}_1$,
$\bar{A}$ be the enumeration of the $\kappa$ Cohen reals.
\begin{myrules}
\item[(1)] In ${\bf V}_1$, there is $(P,\bar{A}) \in K_g$ such that $\Vdash_P
\mbox{``} \gr < \kappa\mbox{''}$, even
$\Vdash_P
\mbox{``} \gr = \aleph_1\mbox{''}$
\item[(2)] For $(P,\bar{A})$ as in (1), we have that
 in ${\bf V}_1$, $\Vdash_P$
``every hyperarithmetic real is
weakly needed for the reaping relation''.
\end{myrules}
\end{claim}
\proof
(1) By \ref{3.5} we have that $\bar{A}$ is
$(g,\kappa)$-o.k.\ in ${\bf V}_1$.
According to \ref{3.7}, we may choose in ${\bf V}_1$
$\lessdot$-increasing and continuous such that 
$(P_i, \bar{A}) \in K$, $P_{i+1}
=
P_i \ast Q_{\name{D^i}}$, where  $\name{D^i} =
\langle \name{D^i_\eta} \such \eta \in {}^{<\omega}
\omega \rangle$ $D^i_\eta = D^i \in V^{P_i}$ as in \ref{3.7}.
Note that $P = \bigcup_{i < \omega_1} P_i$ 
forces that $\gr = \aleph_1$, because it consecutively adds (``shoots'')
$\aleph_1$ reals through ultrafilters in the intermediate models
${\bf V}_0[G_\alpha]$, $\alpha < \omega_1$.
It is easy to see that these $\aleph_1$ reals are a refining family.

(2) 
Now by part (1) and by \ref{crucial} for any $g$
the proof of (2) follows.

\section{There may be more \wneeded reals than
\needed reals}\label{S4} 

Under $\CH$, or if
$||R|| = 2^{\aleph_0}$, then \needed for $R$ and \wneeded
for $R$  coincide.
In this section, we show that there is some quite simply defined relation 
$R$ and that there is some model of $\zfc$ in which 
there are more \wneeded reals for
 $R$ than \needed reals for $R$.
The idea is to use the forcing model from the 
previous section.

\nothing{
\begin{definition}\label{inS1}
\begin{myrules}
\item[(1)]
${\mathcal R}$ is a weak $R$-cover if $(\forall x \in \dom(R))
\; (\exists y \in {\mathcal R}) xRy$.
\item[(2)]
${\mathcal R}$ is an $R$-cover if it is a weak $R$-cover and
$|{\mathcal R}| = ||R||$.
\end{myrules}
\end{definition}
\begin{definition}\label{alt}
We say that $\eta \in {}^\omega 2$ is strongly needed if 
for every weak $R$-cover ${\mathcal R}$ there is some 
$y \in {\mathcal R}$ such that $\eta \leq_T y$.
\end{definition}
}

\begin{claim}\label{4.1}
(Blass \cite{Blass:oberwolfach})
An equivalent condition for ``$\eta \in 
{}^\omega 2$ is \needed for $R$'' is 
$$(\exists x \in \dom(R))(\forall y \in \rge(R))
(x R y 
\rightarrow \eta \leq_T y).$$
\item[(2)] If $2^{\aleph_0} = \aleph_1$ then ``\needed 
for $R$'' is equivalent to ``\wneeded for $R$''
(and for the usual $R$'s, under MA we have that $||R|| =
2^{\aleph_0}$ and hence any adequate set is 
of minimal cardinality and hence the notions coincide).
\end{claim}

\proof 
Suppose that $\eta$ is \needed for $R$ and that there is no
$x$ as in (1).  Then $(\forall x \in \dom(R)) \; (\exists y \in \rge(R))
(xRy \wedge \eta \not\leq_T y)$. So we can build a $R$-adequate set from
all these $y$'s, that shows that $\eta$ is not  \needed for $R$.
For the other implication: Fix $x$ as in (1). Every
 $R$-adequate set has to contain one $y$ such that $xRy$
and hence $\eta \leq_T y$.
\proofend

If $2^{\aleph_0} = \aleph_1$ then ``\needed 
for $R$'' is equivalent to ``\wneeded for $R$''
(and for the usual $R$'s, under MA we have that $||R|| =
2^{\aleph_0}$ and hence any adequate set is 
of minimal cardinality and hence the notions coincide).
But in general, they do not coincide.

\begin{claim}\label{4.2}
There is a simply defined relation $R$ for which it is
consistent that the notions ``\wneededc'' and ``\neededc'' do not
coincide.
In fact, in the forcing model from the previous section, every 
$R$-\needed real is recursive, and all the hyperarithmetic
(and possibly more) reals  are weakly needed for $R$.
\end{claim}
\proof
Let $R= R_0 \cup R_1$, where $R_0$ is the ordinary 
reaping relation, which we write for functions on ${}^\omega 2
\times {}^\omega 2$:
\begin{equation*}
\begin{split}
\eta R_0 \nu & \Leftrightarrow \eta,\nu \in {}^\omega 2
\; \wedge \; (\exists^\infty n) \nu(n) =1 \; \wedge \;
\eta \restriction \nu^{-1}\{1\} \mbox{ is almost constant.}\\
\eta R_1 \nu & \Leftrightarrow \eta,\nu \in {}^\omega 2 \; \wedge \;
(\exists^\infty n) \eta(n)=1 \;\wedge\; 
(\exists^\infty n) \nu(n)=1 \;\wedge\; \\
&\left\langle 
\frac{|\nu^{-1}\{1\} \cap \eta^{-1}\{1\} \cap n |}{|\eta^{-1}\{1\}
\cap n|} \such n \in \omega
\right\rangle \mbox{ converges to } \frac{1}{2}.\\
& \mbox{In particular, for every large enough $n$, } 
\frac{|\nu^{-1}\{1\} \cap \eta^{-1}\{1\} \cap n |}{|\eta^{-1}\{1\}\cap n|}
\in \left[\frac{1}{4}, \frac{3}{4}\right].
\end{split}
\end{equation*}

We use $V^P$ from the previous section. 
There we have that $P = P_0 \ast \name{Q}$, $P_0$ is the forcing 
adding $\kappa$ Cohen reals, and $\name{\bar{A}} $ is an enumeration of
the names of these Cohen reals, and $Q$ is the iteration described in
\ref{3.8}.
Then in $V^P$ we have that
$||R|| \leq ||R_0|| = \aleph_1$.

\smallskip

We first show that every hyperarithmetic real is 
\wneeded for $R$ in this model. 
We take  some $R$-adequate set in $V^P$
 ${\mathcal R}$ of
power 
$\aleph_1$. We let
$$Y_\ell = \{ i < \kappa \such (\exists x \in {\mathcal R}) (A_i R_\ell
x ) \}.$$
So, by the definition of adequate we have that $Y_0 \cup Y_1 = \kappa$.
If $|Y_0| = \kappa$, by the proof of \ref{crucial}, 
we get some $x \in {\mathcal R}$ whose enumeration 
$f$ with $f(n)=m$ if $m$ is the $n$th element of $x$  is
so large in the eventual domination order
that
hyperarithmetic real is computable from it.

We now show that $|Y_1| < \kappa$.
Then it follows that $|Y_0|=\kappa$. Towards a contradiction, we
assume that $|Y_1| = \kappa$. 
In the model from the previous section we 
have that $P= \bigcup_{i < \omega_1} P_i$, $P_0$ adds 
$\kappa$ Cohen reals, $P_i$ increasing and continuous, 
$P_{i+1} = P_i \ast Q_{\name{D_i}}$ as there, $P=
P_0 \ast \name{Q}$. We work in $V^{P_0}$. We have that for some
$p^\ast \in Q/{P_0}$ and some $Q/P_0$-names  
$\name{\nu_i}$, $i < \omega_1$
\[ 
p^\ast \Vdash_{Q/P_0} |\name{Y_1}| = \kappa \; \wedge \;
\name{{\mathcal R}} = \{ \name{\nu_i} \such i < \omega_1 \}.
\]
$Y^\ast = \{ \alpha \such \exists p_\alpha \geq
 p^\ast, p_\alpha \Vdash_{Q/P_0}\alpha \in \name{Y_1}\}.$
By the ccc of $Q/P_0$, we have that $Y^\ast \in [\kappa]^\kappa$,
and for $\alpha \in Y^\ast$ we choose 
$ p^\ast \leq p_\alpha \Vdash_{Q/P_0} 
\mbox{``} \alpha \in \name{Y_1} \mbox{''}$
So for $\alpha \in \kappa$ we have that 
$A_\alpha R_1 \nu_{i(\alpha)}$ and hence for a large 
enough $n^\ast$ for $\kappa$ many 
$\alpha\in Y^\ast$ (w.l.o.g.: for all $\alpha \in Y^\ast$) 
we have that 
$n_\alpha 
= n^\ast$, and there is a $\Delta$-system
for the $\dom{p_\alpha} \in [\kappa \setminus \{0\}]^{<\omega}$ whose root is 
$u^\ast$, $i(\alpha) = i^\ast$.

So we may assume that for $j \in u^\ast$
we have that $p_\alpha(j)$ is an object with trunk $\rho_j$ and not
just a $P_0$-name.
By pure decidability for some $\nu^\ast \in V^{P_0}$ we have:
For every $\alpha \in Y^\ast$ and $m$ 
for some pure extension $q$ of $p_\alpha$
with the same domain  $q \Vdash \name{\nu_{i(\ast)}} \restriction 
m =\nu^\ast \restriction m$. By 
``$n_\alpha = n(\ast)$'' for $\alpha \in Y^\ast$ 
we get an easy contradiction: Suppose $p \in P_0$ and 
\begin{multline*}
p \Vdash_{P_0} \mbox{``}
\forall \alpha \in Y^\ast \;\forall n \geq n(\ast)\;
\exists q_\alpha \geq_{tr} p_\alpha,
\\
q_\alpha \Vdash_{Q/P_0} \mbox{``}\frac{|\name{\nu^\ast}^{-1}\{1\} \cap
A_\alpha^{-1}\{1\} \cap n(\ast)|}{|A_\alpha^{-1}\{1\} \cap n(\ast)|} 
\in \left[\frac{1}{4},\frac{3}{4}\right]\mbox{'' ''.}
\end{multline*}
This is impossible, because we may assume that
$\nu^\ast \in V$ (it needs only countably many of the $\kappa$ Cohen reals)
and we may arrange all other $A_\alpha$'s so that the quotient will be arbitrary.
The forcing $P/P_0$ does not change this fact.

\nothing{
\smallskip
Now we show that if a real is not hyperarithmetic, then it is not
weakly  needed for $R$.
Show: {\em If $\eta$ is weakly needed for $R$, then it is $D$-needed.}
}

\smallskip

Now we show that if a real is not recursive then it is not 
needed for $R$. If $\eta$ is not recursive and 
$x \in {}^\omega 2$, let $\{x,\eta\} \in N \prec (H(\chi),\in)$, $N$
countable. Let $\nu=\nu(x,\eta)$ be random over $N$, and 
we claim 
\begin{equation}\label{claim}
\eta 
\not\leq_{Turing} \nu. 
\end{equation}
Proof of \eqref{claim}: Otherwise we would have
that $\eta$ is recursive in the ground model
by the following:
Suppose
\begin{equation}
\label{computt}
p  \Vdash_{\text{Random}}
\mbox{``$M$ computes $\eta$ from the oracle $\name{\nu}$''}.
\end{equation}

Then by the Lebesgue density theorem we find $s \in
{}^{<\omega} 2$ such that above $s$, $p$ has
Lebesgue measure $> \frac{99}{100} \cdot \Leb(\{\rho
\such s \triangleleft \rho \}$.
The we set
$$B_n = \{ \nu' \in {}^\omega 2 \such s 
\triangleleft \nu' \mbox{ and from $\nu'$ $M$ computes
$\eta(n)$ correctly} \}.$$
From \eqref{computt} we get that $\Leb(B_n) \geq \frac{99}{100}
\cdot \Leb(\{\rho \such s \triangleleft \rho \}$.
So for every sufficiently large $m\in \omega$ we have that
\begin{equation}
\label{decide}
2^{m -\lg(s)} \leq
| \{ \nu' \in {}^m 2 \such s 
\triangleleft \nu' \mbox{ and from $\nu'$ $M$ computes
$\eta(n)$ correctly} \}|.
\end{equation}
So we can run a machine, that has $s$ as an fixed ingredient,
and which, given input $n$,
increases $m$ successively, and then computes $\eta(n)$ with all possible 
oracles above $s$ of length $m\geq\lg(s)$ and decides with
\eqref{decide}, when it is true for  $m$ (and hence for all later $m$), 
which is the right value. So \eqref{claim} is proved.  

But we have that
$x R_1 \nu$ and hence $x R
\nu$. Thus the  collection $\{\nu(x,\eta) \such x \in {}^\omega 2 \}$
is an $R$-adequate family. So there is some $\nu$ such that
$\eta \leq_T \nu$
in contradiction to the equation \eqref{claim}.
So finally we showed that all needed reals for $R$ are recursive.
\proofend

\section{\Needed reals for reaping}\label{S5}

In this section we prove in $\zfc$ that not all hyperarithmetic reals
are \needed for the reaping relation. Since in the model
 from Section 3 all hyperarithmetic
reals are \wneeded for the reaping relation, this model
shows that also for the reaping relation
it is consistent that \wneeded and \needed do not coincide.

\begin{hypothesis}\label{hyp}
We fix $B^\ast \subseteq \omega$ and some $\eta \in {}^\omega 2$ 
such that: if 
$X \subseteq B_\ast = B_1^\ast $ or
$X \subseteq \omega \setminus B_\ast = B_2^\ast$
then $\eta $ is recursive in $\charak_X$.
\end{hypothesis}

By \ref{4.1}, the hypothesis says, that $\eta$ is  
\needed for the reaping relation, with witness $B^\ast$.
For all $X$, that refine $B^\ast$, we have that $\eta$ is recursive in $X$.
Note that \ref{hyp} is similar to $\eta$ being hyperarithmetic: 
the difference is that $\eta$ is computable also 
in every infinite subset of the
complement of $B_\ast$.

\begin{choice}
Let $\langle (M^n_1, M^n_2, a_1^n, a_2^n) \such n < \omega \rangle$ 
be a recursive list of the quadruples $(M_1, M_2, a_1, a_2)$ such that
\begin{myrules}
\item[(i)] $M_1, M_2$ are Turing machines (with reference to an oracle),
\item[(2)] $a_1, a_2$ are finite disjoint sets.
\end{myrules}
W.l.o.g.\ $a_1^n \cup a_2^n \subseteq n$ and each quadruple appears 
infinitely often.
\end{choice}

\begin{definition}\label{5.3}
\begin{myrules}
\item[(1)] We say $\bar{E}=\langle E_n \such n \in \omega \rangle $ is special
if 
\begin{myrules}
\item[(i)] $E_n$ is an equivalence relation on $\omega \setminus n$, and
\item[(ii)] for $m < n$, $E_n$ refines $E_m \restriction
(\omega \setminus n)$,

\item[(iii)] if $A$ is an $E_n$-equivalence class, then 
$A \setminus (n+1)$ is devided by $E_{n+1}$ in at most 
two equivalence classes, and $E_0$ has finitely many classes,

\item[(iv)]    if 

\begin{myrules}
\item[($\alpha$)]
$A$ is an $E_n$-equivalence class and
\item[($\beta$)] there is a partition $X_1, X_2$ of
$A \setminus (n+1)$ such that for all  $j < \omega$, $Y_i \subseteq \omega$, 
$i=1,2$, (if $a_i^n \subseteq Y_i \subseteq X_i \cup a_i^n$,
$h_i < \omega$, the machine $M_i^n$ running with input $j$ and
oracle
$\charak_{Y_i}$ finishes its run giving $h_i$, then 
$h_1 = h_2$),
\end{myrules} 
then $E_{n+1}$ induces such a partition of $A$.
\end{myrules} 
\item[(2)] $\bar{E}$ is special to $\eta$ if in addition

\begin{myrules}
\item[(v)] for all  $A$ and $n$, if $A$ is an $E_n$-class, then
$\eta$ is not recursive 
in $\charak_A$.
\end{myrules} 

\end{myrules} 
\end{definition}

\begin{theorem}\label{5.4}
There is no $\bar{E}$ that is special to $\eta$.
\end{theorem}

\proof
We assume the contrary, and by (Cohen) forcing and
 absoluteness we will derive a contradiction.  The proof will be finished
with \ref{5.11}.

\begin{definition}\label{5.5}
For a special $\bar{E}$ we define $Q= Q_{\bar{E},B_\ast}$ as the following 
notion of forcing:

\begin{myrules}
\item[(1)] $p \in Q$ has the form
$p= (n,A,b_1,b_2) = (n^p, A^p,b_1^p,b_2^p)$ such that

\begin{myrules}
\item[(i)] $n<\omega$,
\item[(ii)] $A$ is an $E_n$-equivalence class,
\item[(iii)] $A$ is infinite,
\item[(iv)] $b_1,b_2$ are disjoint subsets of $n$,
\item[(v)] $b_1 \subseteq B^\ast$, $b_2 \subseteq \omega \setminus B^\ast$.
\end{myrules}

\item[(2)] $p \leq q$ iff

\begin{myrules}
\item[(i)]
 $n^p \leq n^q$, $A^p \supseteq A^q$, $b_i^p \subseteq b_i^q$, for $i=1,2$,
\item[(ii)] $(b^q_1 \cup b_2^q) \setminus (b_1^p \cup b_2^p)
\subseteq A^p$.
\end{myrules}

\item[(3)]
$\name{B_i} = \bigcup \{b^p_i \such p \in 
\name{G_Q} \}$ is a $Q$-name of a subset $B_i \in V[G]$
 of $B_i^\ast$ if $i=1,2$.
\end{myrules}
\end{definition}

So  if $E_0$ has finitely many equivalence classes, then
$Q$ is equivalent to Cohen forcing and independent of
$\bar{E}$ and $B_\ast$. Nevertheless we keep 
the complicated conditions, because they are better tailored
for $\eta$'s needed for the reaping relation. 

\begin{claim}\label{5.6}
For $i=1,2$ we have
\begin{myrules}
\item[(1)] $\Vdash_Q \mbox{``}\name{b_i}$  is an infinite
subset of $B^\ast_i$''.
\item[(2)] For some $p^\ast$, $p^\ast \Vdash_Q \mbox{``}M^{n^{p^\ast}}_i$ 
computes $\eta$ with the oracle $\charak_{\name{b_i}}$''.
\end{myrules}
\end{claim}
\proof
(1) Fix $i=1$ or $i=2$. 
It is enough, to find for a given $p \in Q$ some
$q \geq p$, $q \in Q$ such that $b_i^p \subsetneqq b_i^q$.
Now $A^p \cap B_i^\ast$ is infinite, because of
the hypothesis on $B^\ast$ and because
$\eta$ is not recursive in $\charak_{A^p}$
by the assumption of the indirect proof of \ref{5.4}. 
We may choose $h \in A^p \cap B_i^\ast$,
$h \geq \max(b_i^p)+1$ and an 
infinite $E_{h+1}$-class $A \subseteq A^p$, which exists
because $A^p$ is infinite and because $E_{h+1}$ has finitely many
equivalence classes.
We define $q$ as
$n^q= h+1$, $A^q = A$, $b_i^q = b_i^p \cup \{h\}$, $b^q_{3-i} = b^p_{3-i}$.

(2) The statement made in Hypothesis~\ref{hyp} 
on $B^\ast$ and on $\eta$ is $\Pi_1^1$ and holds in $V$, hence it
holds in $V[G]$ as well by \cite[Theorem 98, p.\ 530]{Jech}. 
Now we apply it in $V[G]$ to part (1) of this claim.
\proofend

\smallskip

We fix $p^\ast$, $M_1^{n^{p^\ast}}$,
 $M_2^{n^{p^\ast}}$ as  in part (2) of Claim \ref{5.6}.

\begin{fact}\label{5.7} 
There is some $q \geq p^\ast$ such that
for $i=1,2$, $M_i^{n^q} = M_i^{n^{p^\ast}}$
and such that $b_i^q=a_i^{n^q}$.
\end{fact}

\proof For some $n^\ast \geq n^{p^\ast}$ the quadruple
$(M_1^{n^\ast}, M_2^{n\ast},
a_1^{n^\ast}, a_2^{n^\ast})$ is equal to\\
$(M_1^{n^{p^\ast}}, M_2^{n^{p^\ast}}, b_1^{p^\ast}, b_2^{p^\ast})$. 
Let $A$ be an infinite
$E_{n^\ast}$-class which is a subset of $A^{p^\ast}$.
So we take $q=
(n^\ast, A, a_1^{n^\ast}, a_2^{n^\ast})$.

\begin{claim}\label{5.8}
For $n^\ast$, $A$ the demands $(\alpha) + (\beta)$ of clause 
(iv) of \ref{5.3} hold, hence the conclusion.
\end{claim}

\proof
We work first in $V[G]$. There, by \ref{5.6},
$X_i = A \cap B_i^\ast$ and $A$ exemplify \ref{5.3}(iv).
But \ref{5.3}(iv) is a $\Sigma_2^1$-statement of the parameters 
$(A, a_1^n,a_2^n)$, and therefore it holds in $V$ as well by
Shoenfield's absoluteness theorem \cite[Theorem 98, p.\ 530]{Jech}.
\proofend

\begin{convention}\label{5.9}
Let $A_1 \neq A_2$ be the $E_{n^\ast+1}$-equivalence classes
which are subsets of $A$, with $A_i$ for $M_i$
as in \ref{5.3}(iv).
\end{convention}

\begin{claim}\label{5.10}
If $j< \omega$ then for some $b \subseteq m < \omega$ we have that 
$b \cap n^{q} = b_1^{q}$, $b \setminus n^{q} \subseteq A_1$,
if we let $M_1$ run with input $j$ and oracle $\charak_b \restriction m$ 
it gives an answer 
(i.e.\ it finishes and asks the oracle only questions
in its domain $m$).
\end{claim}

\proof:
We define $r \in Q$ by $n^r = n^{q} +1$, $A^r = A_1$, $b_i^r =
b_i^q$ for $i = 1,2$.
So $q \leq r \in Q$.  By the choice of $q$ for some $s \in Q$, 
$r \leq s$ and $s$ forces a value to the run of $M_1$ with input $j$ 
and oracle $\name{b_1}$ so also to the answers to the oracle in this run. 
Let $b = b_1^s$.
\proofend

\begin{claim}\label{5.11}
For every $j \in \omega$, $k \in 2$ the following are equivalent
\begin{myrules}
\item[(1)] $\eta(j)=k$.
\item[(2)]  For some $b \subseteq m < \omega$, $b\cap n^{q}
= b_1^{q} $, $b \setminus n^{p^\ast} \subseteq A_1$, and 
$M_1$ running with input $j$ and oracle $\charak_b \restriction m$ 
gives the answer $k$.
\end{myrules}
\end{claim}

\proof
$(i) \rightarrow (ii)$ by \ref{5.10}.
Since $\eta \in V$, the reverse implication holds as well.
\nothing{
If $(ii) + \neg (i)$, then we apply (ii) to get $b_1^k$ and 
\ref{5.9} and $\eta(j) \neq k$ to get $b_1^{3-k}$. 
This together with the choice of 
$E_{n^{q} +1 }$ and the Definition~\ref{5.3}(iv)($\beta$) gives
a contradiction.}
\proofend

\smallskip

End of the proof of \ref{5.4}:
$\eta$ is recursive in $\charak_{A_1}$.
By \ref{5.11} we  try all $b$'s for a given $j$
and hence $\eta$ is recursive in $\charak_{A_1}$.
How to run through all trials is explained in more detail
in \cite[Theorem 9]{bl-ober}.
\proofendof{\ref{5.4}}

\begin{claim}\label{5.12}
There is a special $\bar{E}$ that has as a three place relation 
$\{ \langle n,x,y \rangle \such
x E_n y \}$
Turing degree $\leq O^\omega$ and such that 
if $A$ is an $E_n$-equivalence  
class then $\charak_A \leq_T O^{n+1}$, the $(n+1)$st jump of $O$.
\end{claim}

\proof We choose $E_n$ by induction on $n$. 

$n=0$. If for every $m$ there is a partition $(c_0,c_1)$ of $m$ such that
for $i \in \{1,2\}$ for every $b_i \subseteq c_i$ and $j <n$
and $M^0_i$ running with input $j$ and oracle $\charak_{b_i}\restriction m$
or $\charak_{b_i}\restriction m$ and giving the results
$k_i$ then $k_0 = k_1$, then 
we choose among these pairs $(c_1^m,c_2^m)$ such that
$\charak_{c_1^m}$ is minimal in the lexicographical order.
If $(c_1^m, c_2^m)$ are defined for every $m$, then we have that
$m^1 \leq m^2 \leq m^3 \Rightarrow 
\charak_{c_1^{m_2} \cap m_1 }
\leq_{lex}
\charak_{c_1^{m_3}\cap m_1}$.
So $\langle c_1^m \such m \in \omega \rangle$ converges to some $c_1$.
Now we define $E_0$, having two classes:
$c_1$ and $\omega \setminus c_1$.
The relation $E_0$ is computable in $O^1$.

In the step from $n$ to $n+1$,  the relation
$E_{n+1}$ is defined similarly, with the modification that
we use the description of $E_n$ as a parameter and take partitions
$(c_0,c_1)$ of $(m\setminus n) \cap C$ for each $E_n$-class $C$, and
oracles $B_i \cup a^n_i$.
Clearly using $O^{n+2}$ we can choose $E_{n+1}$ and $\bar{E}$ is $\ll \Delta_1^1$.
\proofend

\begin{remark*}
Just to show that Con(\needed for reaping 
does not coincide with \wneeded for reaping) is is enough to find a 
$\Delta_1^1$-relation $\bar{E}$ which is special.
\end{remark*}

\begin{conclusion}\label{5.13}
If $\eta$ is \needed for the reaping relation, then $\bigvee_{n \in \omega} 
(\eta \leq_T O^n)$, hence in the $V^P$ from Section~\ref{S3}
 many $\Delta_1^1$ reals are 
not \needed for the reaping relation, but only \wneeded for the
reaping relation.
\end{conclusion}

\proof We take $\bar{E}$ as in \ref{5.12}. From \ref{5.11} we get that
$\bar{E}$ is not special to any $\eta$ 
that satisfies $\ref{hyp}$ for some $B_\ast$. 
So any $\eta$ that is \needed for the reaping relation is recursive in 
$\bar{E}$.
\proofend 

\section{Coincidence} \label{S6}
\nc{\random}{{\text{random}}}
In this section we give a condition on a relation $R$
 under which 
\needed for $R$ and \wneeded for $R$ coincide and 
show that the condition is fulfilled for
the relation $R$ defined below.

\begin{definition}\label{6.1}
The domain of the relation $R_{\random}$ is 
$\{ \eta \such \eta$ is a code for a measure 1 set, say a tree
$T_\eta \subseteq {}^{<\omega} 2$ of positive measure$\}$. 
The range of $R_\random$ is ${}^\omega 2$. We set
$\eta R_\random \nu$ iff $\nu \in A_\eta :=
\{\rho \in {}^\omega 2 
\such $ for some $\rho' \in T_\eta$ we have that $\rho=^\ast \rho' \}$.
\end{definition}

\begin{claim}\label{6.2}
 \begin{myrules}
\item[(1)] Assume that 
\begin{equation}\tag{$\otimes_R$}
\begin{split}
&\mbox{(a) $R$ is a 2-place Borel relation on 
${}^\omega 2$, and }\\
&\mbox{(b) for every $x_1,x_2 \in {}^\omega 2 $,
if $x_2$ is not recursive, there is $x \in {}^\omega 2$ such that }\\
& \makebox[2cm]{} \otimes \makebox[1cm]{}
(\forall \nu) \biggl(x R \nu \rightarrow (x_1 R \nu \wedge \neg (x_2 
\leq_T \nu))\biggr).
\end{split}
\end{equation}
then the notions of strongly needed for $R$ and weakly needed
for $R$ coincide
and coincide with being recursive.
\item[(2)]
 The relation $R_\random$ satisfies the criterion $\otimes_R$ from
Part (1).
\end{myrules}
\end{claim}

\proof
(1) We have show that every weakly needed real is recursive.
Then by ``recursive $\rightarrow$ strongly needed $\rightarrow$
weakly needed $\rightarrow$ recursive'' all three notions coincide.

Suppose that $x^\ast \in {}^\omega 2$ is not recursive.
We show that $x^\ast$ is not weakly needed. 
Let $Y$ be a strong $R$-cover. Let $Y^\ast =
\{ \nu \in Y \such 
\neg x^\ast \leq_T \nu \}$. $Y^\ast \subseteq Y$, and hence
$|Y^\ast| \leq |Y| = ||R||$. We show that $Y^\ast$ is also an $R$-cover.
Let $x_1 \in {}^\omega 2$ be given. 
We take $x_2 = x^\ast$, and apply (b) of $\otimes_R$. So we get 
$x$ as there. Since $Y$ is an $R$-cover we find some 
$\nu \in Y$ such that $x R \nu$.
Hence by $\otimes_R$ we have that $x_1 R \nu \wedge x_2 \not\leq_T \nu$.
So $\nu \in Y^\ast$ $R$-covers $x_1$. 

(2) Let $x_1, x_2$ be given. We take $N \prec (H(\beth_3),\in)$ such that
$x_1,x_2 \in N$. Let $T = T_\eta$ be Amoeba-generic over $N$.
Then $\eta= x$ is as claimed in (1)(b).

\proofend

\begin{conclusion}
Strongly $R_\random$-needed and weakly $R_\random$-needed coincide
and are just all the recursive reals. \proofend
\end{conclusion}


\begin{thebibliography}{1}

\bibitem{abraham:handbook}
Uri Abraham.
\newblock Proper forcing.
\newblock In Matthew Foreman, Akihiro Kanamori, and Menachem Magidor, editors,
  {\em Handbook of Set Theory}. Kluwer, To appear.

\bibitem{BJ}
Tomek Bartoszy\'{n}ski and Haim Judah.
\newblock {\em {{Set Theory, On the Structure of the Real Line}}}.
\newblock A K Peters, Wellesley, Massachusetts, 1995.

\bibitem{Blass:oberwolfach}
Andreas Blass.
\newblock Needed reals.
\newblock {\em Talk at Oberwolfach December 1999}.

\bibitem{bl-ober}
Andreas Blass.
\newblock Needed reals and recursion in generic reals.
\newblock {\em To appear in APAL},
  2001.

\bibitem{Jech}
Thomas Jech.
\newblock {\em Set Theory}.
\newblock Addison Wesley, 1978.

\bibitem{Jockusch}
Carl G.~Jockusch Jr.
\newblock Uniformly introreducible sets.
\newblock {\em J. Symbolic Logic}, 33:521--536, 1968.

\bibitem{Oxtoby}
John Oxtoby.
\newblock {\em {Measure and Category}}.
\newblock Springer, second edition, 1980.

\bibitem{Sh:707}
Saharon Shelah.
\newblock Tree forcings.
\newblock {\em Preprint [Sh:707]}, 2000.

\bibitem{solovay:hyper}
Robert Solovay.
\newblock Hyperarithmetically computable sets.
\newblock {\em Trans. Amer. Math. Soc.}, 239:99--122, 1978.

\end{thebibliography}
%
\def\germ{\frak} \def\scr{\cal}
  \ifx\documentclass\undefinedcs\def\rm{\fam0\tenrm}\fi
  \def\defaultdefine#1#2{\expandafter\ifx\csname#1\endcsname\relax
  \expandafter\def\csname#1\endcsname{#2}\fi} \defaultdefine{Bbb}{\bf}
  \defaultdefine{frak}{\bf} \defaultdefine{mathbb}{\bf}
  \defaultdefine{beth}{BETH} \def\bbfI{{\Bbb I}} \def\mbox{\hbox}
  \def\text{\hbox} \def\om{\omega} \def\Cal#1{{\bf #1}} \def\pcf{pcf}
  \defaultdefine{cf}{cf} \defaultdefine{reals}{{\Bbb R}}
  \defaultdefine{real}{{\Bbb R}} \def\restriction{{|}} \def\club{CLUB}
  \def\w{\omega} \def\exist{\exists} \def\se{{\germ se}} \def\bb{{\bf b}}
  \def\equivalence{\equiv} \let\lt< \let\gt> \def\cite#1{[#1]}

\end{document}